\documentclass[12pt]{amsart}

\usepackage{latexsym,pstricks}

\newcommand\courbes{\psset{unit=.8cm}
\begin{pspicture}(-4.2,-5)(3.9,4)
\psset{linewidth=1pt}

\psarc(0,1.2){1.5}{230}{-50}
\psarc(0,-.8){1}{45}{-225}

\pscustom{\rotate{-5}
\pscurve(3,-1.8)(3.2,-3.2)(1.8,-3)}
\pscustom{\rotate{-5}
\pscurve(3.28,-1.2)(3.6,-3.6)(1.2,-3.28)}
\pscustom{\rotate{-5}
\pscurve(3.28,-1.2)(3.15,-1.15)(3.05,-1.2)(2.9,-1.4)
(2.87,-1.5)(2.9,-1.67)(3,-1.8)}
\pscustom[linestyle=dotted,dotsep=1pt]{\rotate{-5}
\pscurve(3.28,-1.2)(3.15,-1.3)(3.05,-1.5)(3,-1.8)}
\pscustom{\rotate{-43.8} \pscurve(3.28,-1.2)(3.15,-1.15)
(3.05,-1.2)(2.9,-1.4)(2.87,-1.5)(2.9,-1.67)(3,-1.8)}
\pscustom[linestyle=dotted,dotsep=1pt]{\rotate{-43.8}
\pscurve(3.28,-1.2)(3.15,-1.3)(3.05,-1.5)(3,-1.8)}

\pscustom{\rotate{-80}
\pscurve(3,-1.8)(3.2,-3.2)(1.8,-3)}
\pscustom{\rotate{-80}
\pscurve(3.28,-1.2)(3.6,-3.6)(1.2,-3.28)}
\pscustom{\rotate{-80}
\pscurve(3.28,-1.2)(3.15,-1.15)(3.05,-1.2)(2.9,-1.4)
(2.87,-1.5)(2.9,-1.67)(3,-1.8)}
\pscustom[linestyle=dotted,dotsep=1pt]{\rotate{-80}
\pscurve(3.28,-1.2)(3.15,-1.3)(3.05,-1.5)(3,-1.8)}
\pscustom{\rotate{-118.8} \pscurve(3.28,-1.2)(3.15,-1.15)
(3.05,-1.2)(2.9,-1.4)(2.87,-1.5)(2.9,-1.67)(3,-1.8)}
\pscustom[linestyle=dotted,dotsep=1pt]{\rotate{-118.8}
\pscurve(3.28,-1.2)(3.15,-1.3)(3.05,-1.5)(3,-1.8)}

\pscustom{\rotate{-175}
\pscurve(3,-1.8)(3.2,-3.2)(1.8,-3)}
\pscustom{\rotate{-175}
\pscurve(3.28,-1.2)(3.6,-3.6)(1.2,-3.28)}
\pscustom{\rotate{-175}
\pscurve(3.28,-1.2)(3.15,-1.15)(3.05,-1.2)(2.9,-1.4)
(2.87,-1.5)(2.9,-1.67)(3,-1.8)}
\pscustom[linestyle=dotted,dotsep=1pt]{\rotate{-175}
\pscurve(3.28,-1.2)(3.15,-1.3)(3.05,-1.5)(3,-1.8)}
\pscustom{\rotate{-213.8} \pscurve(3.28,-1.2)(3.15,-1.15)
(3.05,-1.2)(2.9,-1.4)(2.87,-1.5)(2.9,-1.67)(3,-1.8)}
\pscustom[linestyle=dotted,dotsep=1pt]{\rotate{-213.8}
\pscurve(3.28,-1.2)(3.15,-1.3)(3.05,-1.5)(3,-1.8)}

\psarc(3.7,0){.45}{120}{240}
\psarc(4.65,0){1.25}{160}{200}

\pscustom{\rotate{30}
\psarc(3.7,0){.45}{120}{240}}
\pscustom{\rotate{30}
\psarc(4.65,0){1.25}{160}{200}}

\pscustom{\rotate{90}
\psarc(3.7,0){.45}{120}{240}}
\pscustom{\rotate{90}
\psarc(4.65,0){1.25}{160}{200}}

\psarc(0,0){3.5}{6.5}{23}
\psarc(0,0){3.5}{37}{47}
\psarc[linestyle=dotted,dotsep=3pt](0,0){3.5}{48}{72}
\psarc(0,0){3.5}{73}{83}
\psarc(0,0){3.5}{97}{116}
\psarc(0,0){3.5}{126}{154}
\psarc(0,0){3.5}{165}{175}
\psarc[linestyle=dotted,dotsep=3pt](0,0){3.5}{176}{199}
\psarc(0,0){3.5}{200}{210}
\psarc(0,0){3.5}{221}{249}
\psarc(0,0){3.5}{260}{285}
\psarc(0,0){3.5}{296}{324}
\psarc(0,0){3.5}{335}{353.5}

\psset{linewidth=.5pt}

\psellipse(0,0)(1.6,1)
\rput[l]{0}(-1.55,0){$\beta$}

\pscustom{\rotate{-5}
\pscurve(3.2,-1.6)(3.4,-3.4)(1.6,-3.2)}
\pscustom{\rotate{-5}
\pscurve(1.6,-3.2)(1.3,-2.9)(1.2,-2.5)(1.5,-1.5)(2.4,-1.2)(3,-1.3)(3.2,-1.6)}
\rput[l]{0}(1,-2.6){$\beta_{_{1}}$}

\pscustom{\rotate{-80}
\pscurve(3.2,-1.6)(3.4,-3.4)(1.6,-3.2)}
\pscustom{\rotate{-80}
\pscurve(1.6,-3.2)(1.3,-2.9)(1.2,-2.5)(1.5,-1.5)(2.4,-1.2)(3,-1.3)(3.2,-1.6)}
\rput[l]{0}(-2.5,-1.95){$\beta_{_{2}}$}

\pscustom{\rotate{-175}
\pscurve(3.2,-1.6)(3.4,-3.4)(1.6,-3.2)}
\pscustom{\rotate{-175}
\pscurve(1.6,-3.2)(1.3,-2.9)(1.2,-2.5)(1.5,-1.5)(2.4,-1.2)(3,-1.3)(3.2,-1.6)}
\rput[l]{0}(-3,.65){$\beta_{_{g-1}}$}

\pscurve(.7,-.1)(1.6,-.2)(2,-.3)(2.5,-.45)(2.9,-.6)(3.3,-.9)(3.35,-1)
\pscurve[linestyle=dashed](.7,-.1)(.75,-.2)(1.2,-.5)
(2,-.85)(2.5,-1)(3,-1.1)(3.35,-1)
\rput[lb]{0}(2.2,-.9){$\alpha_{_{1}}$}

\pscurve(2.2,-2.71)(2.23,-2.6)(2.15,-2.3)(2,-2)(1.6,-1.5)
(1.2,-1)(.9,-.6)(.42,-.22)
\pscurve[linestyle=dashed](.42,-.22)(.45,-.4)(.5,-.6)(.67,-1)(1,-1.55)
(1.4,-2.1)(1.7,-2.45)(1.9,-2.6)(2.2,-2.71)
\rput[lb]{0}(1.42,-2.18){$\alpha_{_{2}}$}

\pscurve(.1,-3.5)(.22,-3.3)(.3,-3)(.4,-2.6)(.5,-1.8)
(.45,-1.3)(.4,-1)(.3,-.6)(.2,-.29)
\pscurve[linestyle=dashed](.2,-.29)(.1,-.45)(0,-.7)
(-.1,-1)(-.2,-2)(-.15,-3.1)(0,-3.4)(.1,-3.5)
\rput[rb]{0}(.43,-1.9){$\alpha_{_{3}}$}

\pscurve(-1.95,-2.9)(-1.65,-2.8)(-1.3,-2.45)(-1,-2)
(-.7,-1.55)(-.4,-.9)(-.2,-.3)
\pscurve[linestyle=dashed](-.2,-.3)(-.4,-.4)(-1,-1)
(-1.4,-1.5)(-1.7,-1.9)(-1.9,-2.3)(-1.97,-2.7)(-1.95,-2.9)
\rput[rb]{0}(-.7,-1.6){$\alpha_{_{4}}$}

\pscurve(-3.16,-1.5)(-2.9,-1.55)(-2.6,-1.5)(-2.3,-1.4)(-1.3,-.9)
(-.6,-.4)(-.5,-.2)
\pscurve[linestyle=dashed](-.5,-.2)(-.9,-.25)(-1.2,-.3)(-1.7,-.5)
(-2.2,-.7)(-2.7,-1)(-3,-1.2)(-3.1,-1.4)(-3.16,-1.5)
\rput[rb]{0}(-2.3,-1.4){$\alpha_{_{5}}$}

\pscurve(-2.61,2.31)(-2.58,2.1)(-2.4,1.8)(-2,1.4)(-1,.5)(-.45,.1)
\pscurve[linestyle=dashed](-.45,.1)(-.5,.3)(-.58,.5)(-.7,.7)(-.9,1)(-1.9,2)
(-2.2,2.25)(-2.5,2.32)(-2.61,2.31)
\rput[lb]{0}(-2.49,1.78){$\alpha_{_{2g-2}}$}

\pscurve(-.92,3.38)(-1,3.2)(-1.05,3)(-1,2.6)(-.9,2)(-.7,1.4)(-.5,.8)(-.2,.2)
\pscurve[linestyle=dashed](-.2,.2)(-.15,.6)(-.2,1.3)(-.3,2)(-.4,2.4)(-.6,3)
(-.7,3.2)(-.92,3.38)
\rput[l]{0}(-.95,2.6){$\alpha_{_{2g-1}}$}

\pscurve(.7,3.42)(.5,3.25)(.4,3.05)(.3,2.8)(.1,1.9)(.05,1.4)(.03,.7)
(.1,.5)(.2,.19)
\pscurve[linestyle=dashed](.2,.19)(.3,.4)(.6,1)(.8,1.7)(.85,2)(.89,2.5)
(.85,3)(.8,3.2)(.7,3.42)
\rput[lb]{0}(.2,2.2){$\alpha_{_{2g}}$}

\pscurve(3.35,1)(3.2,1.05)(2.9,1.1)(2.4,1)(1.7,.8)(1,.5)(.6,.25)(.45,.1)
\pscurve[linestyle=dashed](.45,.1)(.7,.05)(1,0)(1.6,.1)(2.3,.3)(2.8,.5)
(3.1,.7)(3.3,.9)(3.35,1)
\rput[lt]{0}(1.8,.87){$\alpha_{_{2g+n-2}}$}

\end{pspicture}
\begin{pspicture}(-3.8,-5)(3.7,4)

\psset{linewidth=1pt}

\psarc(0,1.2){1.5}{230}{-50}
\psarc(0,-.8){1}{45}{-225}

\pscustom{\rotate{40}
\psarc(3.7,0){.45}{120}{240}}
\pscustom{\rotate{40}
\psarc(4.65,0){1.25}{160}{200}}

\pscustom{\rotate{160}
\psarc(3.7,0){.45}{120}{240}}
\pscustom{\rotate{160}
\psarc(4.65,0){1.25}{160}{200}}

\pscustom{\rotate{280}
\psarc(3.7,0){.45}{120}{240}}
\pscustom{\rotate{280}
\psarc(4.65,0){1.25}{160}{200}}

\psarc(0,0){3.5}{46.4}{65}
\psarc[linestyle=dotted,dotsep=3pt](0,0){3.5}{69}{129}
\psarc(0,0){3.5}{133}{153}
\psarc(0,0){3.5}{166.5}{186}
\psarc[linestyle=dotted,dotsep=3pt](0,0){3.5}{190}{250}
\psarc(0,0){3.5}{254}{273.5}
\psarc(0,0){3.5}{286.5}{306}
\psarc[linestyle=dotted,dotsep=3pt](0,0){3.5}{310}{10}
\psarc(0,0){3.5}{14}{33.6}

\psset{linewidth=.5pt}

\pscurve(1.56,-3.12)(1.6,-2.9)(1.55,-2.6)(1.4,-2.2)(1,-1.4)(.5,-.7)(.1,-.3)
\pscurve[linestyle=dashed](.1,-.3)(.05,-.4)(.02,-.55)(.1,-.8)
(.2,-1.1)(.6,-2)(.9,-2.5)(1.2,-2.9)(1.3,-3)(1.56,-3.12)
\rput[lt]{0}(1.1,2.3){$\alpha_{_{i}}$}

\pscustom{\rotate{-167}
\pscurve(.67,-.09)(1.6,-.2)(2,-.3)(2.5,-.45)(2.9,-.6)(3.3,-.9)(3.35,-1)}
\pscustom[linestyle=dashed]{\rotate{-167}
\pscurve(.67,-.09)(.75,-.2)(1.2,-.5)
(2,-.85)(2.5,-1)(3,-1.1)(3.35,-1)}
\rput[rt]{0}(1,-1.4){$\alpha_{_{j}}$}

\pscurve(1.8,3)(1.4,2.7)(.95,2.05)(.7,1.5)(.2,.2)
\pscurve[linestyle=dashed](.2,.2)(.3,.23)(.55,.4)(.9,.8)(1.2,1.2)(1.5,1.7)
(1.7,2.2)(1.8,2.5)(1.88,2.8)(1.85,2.9)(1.8,3)
\rput[lb]{0}(-2.9,-.05){$\alpha_{_{k}}$}

\pscurve(1.2,-3.28)(1.15,-3.1)(.85,-2.7)(.2,-2)(-.3,-1.6)(-1.5,-.9)(-2.4,-.5)
(-3,-.3)(-3.49,-.2)
\pscurve[linestyle=dashed](-3.49,-.2)(-3,-.7)(-2.5,-1.1)(-1.7,-1.7)(-.7,-2.3)
(.4,-2.9)(1,-3.2)(1.2,-3.28)
\rput[l]{0}(1.6,0){$\gamma_{_{i,j}}$}

\pscustom{\rotate{-120}
\pscurve(1.2,-3.28)(1.15,-3.1)(.85,-2.7)(.2,-2)(-.3,-1.6)(-1.5,-.9)(-2.4,-.5)
(-3,-.3)(-3.49,-.2)}
\pscustom[linestyle=dashed]{\rotate{-120}
\pscurve(-3.49,-.2)(-3,-.7)(-2.5,-1.1)(-1.7,-1.7)(-.7,-2.3)
(.4,-2.9)(1,-3.2)(1.2,-3.28)}
\rput[t]{0}(-.8,-1.4){$\gamma_{_{j,k}}$}

\pscustom{\rotate{120}
\pscurve(1.2,-3.28)(1.15,-3.1)(.85,-2.7)(.2,-2)(-.3,-1.6)(-1.5,-.9)(-2.4,-.5)
(-3,-.3)(-3.49,-.2)}
\pscustom[linestyle=dashed]{\rotate{120}
\pscurve(-3.49,-.2)(-3,-.7)(-2.5,-1.1)(-1.7,-1.7)(-.7,-2.3)
(.4,-2.9)(1,-3.2)(1.2,-3.28)}
\rput[rb]{0}(-.4,1.6){$\gamma_{_{k,i}}$}

\end{pspicture}}

\newcommand\etoile{\psset{unit=.5cm}
\begin{pspicture}(-4,-4)(4,4)

\psset{linewidth=1pt}

\psarc(0,1.2){1.5}{230}{-50}
\psarc(0,-.8){1}{45}{-225}

\psccurve(3.45,.6)(3.3,.5)(3.2,.35)(3.1,0)(3.2,-.35)(3.3,-.5)(3.45,-.6)
(3.4,-.5)(3.35,-.3)(3.3,0)(3.35,.3)(3.4,.5)

\pscustom{\rotate{120}
\psccurve(3.45,.6)(3.3,.5)(3.2,.35)(3.1,0)(3.2,-.35)(3.3,-.5)(3.45,-.6)
(3.4,-.5)(3.35,-.3)(3.3,0)(3.35,.3)(3.4,.5)}

\pscustom{\rotate{-120}
\psccurve(3.45,.6)(3.3,.5)(3.2,.35)(3.1,0)(3.2,-.35)(3.3,-.5)(3.45,-.6)
(3.4,-.5)(3.35,-.3)(3.3,0)(3.35,.3)(3.4,.5)}

\psarc(0,0){3.5}{10}{110}
\psarc(0,0){3.5}{130}{230}
\psarc(0,0){3.5}{250}{-10}

\psset{linewidth=.5pt}

\pscurve(1.56,-3.12)(1.6,-2.9)(1.55,-2.6)(1.4,-2.2)(1,-1.4)(.5,-.7)(.1,-.3)
\pscurve[linestyle=dashed](.1,-.3)(.05,-.4)(.02,-.55)(.1,-.8)
(.2,-1.1)(.6,-2)(.9,-2.5)(1.2,-2.9)(1.3,-3)(1.56,-3.12)
\rput[lt]{0}(1.25,2.5){$\alpha_{_{3}}$}

\pscustom{\rotate{-167}
\pscurve(.67,-.09)(1.6,-.2)(2,-.3)(2.5,-.45)(2.9,-.6)(3.3,-.9)(3.35,-1)}
\pscustom[linestyle=dashed]{\rotate{-167}
\pscurve(.67,-.09)(.75,-.2)(1.2,-.5)
(2,-.85)(2.5,-1)(3,-1.1)(3.35,-1)}
\rput[t]{0}(1.0,-1.65){$\alpha_{_{1}}$}

\pscurve(1.8,3)(1.4,2.7)(.95,2.05)(.7,1.5)(.2,.2)
\pscurve[linestyle=dashed](.2,.2)(.3,.23)(.55,.4)(.9,.8)(1.2,1.2)(1.5,1.7)
(1.7,2.2)(1.8,2.5)(1.88,2.8)(1.85,2.9)(1.8,3)
\rput[lt]{0}(-3.3,-.15){$\alpha_{_{2}}$}

\pscurve(3.2,1.4)(3,1.2)(2.8,1)(2.5,.6)(2.4,.4)(2.3,0)(2.4,-.4)(2.5,-.6)
(2.8,-1)(3,-1.2)(3.2,-1.4)
\pscurve[linestyle=dashed](3.2,-1.4)(3,-1)(2.75,-.3)(2.7,0)
(2.75,.3)(3,1)(3.2,1.4)
\rput[rt]{0}(2.95,-.8){$\gamma_{_{1}}$}

\pscustom{\rotate{-120}
\pscurve(3.2,1.4)(3,1.2)(2.8,1)(2.5,.6)(2.4,.4)(2.3,0)(2.4,-.4)(2.5,-.6)
(2.8,-1)(3,-1.2)(3.2,-1.4)}
\pscustom[linestyle=dashed]{\rotate{-120}
\pscurve(3.2,-1.4)(3,-1)(2.75,-.3)(2.7,0)
(2.75,.3)(3,1)(3.2,1.4)}
\rput[rb]{0}(-2.2,-1.85){$\gamma_{_{2}}$}

\pscustom{\rotate{120}
\pscurve(3.2,1.4)(3,1.2)(2.8,1)(2.5,.6)(2.4,.4)(2.3,0)(2.4,-.4)(2.5,-.6)
(2.8,-1)(3,-1.2)(3.2,-1.4)}
\pscustom[linestyle=dashed]{\rotate{120}
\pscurve(3.2,-1.4)(3,-1)(2.75,-.3)(2.7,0)
(2.75,.3)(3,1)(3.2,1.4)}
\rput[lt]{0}(-2.8,1.9){$\gamma_{_{3}}$}

\psellipse(0,0)(1.3,1)
\rput[rb]{0}(-.7,.75){$\beta$}

\end{pspicture}}

\renewcommand\bar[1]{\overline{#1}}

\renewcommand\proof{\noindent{\bf Proof.\ }\,}
\newcommand\eproof{\hfill $\Box$}

\input diagrams
\newarrow{to} ---->

\newtheorem{theorem}{Theorem}
\newtheorem{proposition}{Proposition}
\newtheorem{lemma}[proposition]{Lemma}
\newtheorem{rem}{Remark}

\newenvironment{remark}{\noindent\begin{rem}{\bf\hspace*{-3mm}} \rm}
                       {\end{rem}\vskip3mm}

\newcounter{liste}

\begin{document}

\title[Presentation of Mapping Class groups]{A finite presentation of
the mapping class group of an oriented surface}

\author{Sylvain GERVAIS}

\address{Universit\'e de Nantes-UMR 6629, 2 rue de la Houssini\`ere,
BP 92208, 44322 NANTES Cedex 3, FRANCE}
\email{gervais$\char'100$math.univ-nantes.fr}
\date{November 23, 1998.}

\begin{abstract}
We give a finite presentation of the mapping class group of an oriented 
(possibly bounded) surface of genus greater or equal than $1$, considering Dehn
twists on a very simple set of curves.
\end{abstract}
\maketitle


\section*{Introduction and notations}

Let $\,\Sigma_{g,n}\,$ be an oriented surface of genus $\,g\!\geq\! 
1\,$ with $n$ boundary components and denote by $\,\mathcal{M}_{g,n}\,$
its mapping class group, that is to say the group of orientation preserving 
diffeomorphisms of $\,\Sigma_{g,n}\,$ which are the identity on 
$\,\partial\Sigma_{g,n}$, modulo isotopy: 
$$\,\mathcal{M}_{g,n}=\pi_{0}\bigl(\hbox{Diff}^{+}(\Sigma_{g,n},
\partial\Sigma_{g,n})\bigr)\,.$$

For a simple closed curve $\alpha$ in $\,\Sigma_{g,n}$, denote by
$\tau_{\alpha}\,$ the Dehn twist along $\alpha$. If $\alpha$ and $\beta$ are 
isotopic, then the associated twists are also isotopic: thus, we shall consider 
curves up to isotopy. We shall use greek letters to denote them, and we 
shall not distinguish a Dehn twist from its isotopy class.

It is known that $\,\mathcal{M}_{g,n}\,$ is generated by Dehn twists
\cite{Dehn,Lickorish1,Lickorish2}.\linebreak[4] Wajnryb 
gave in \cite{Wajnryb} a presentation of 
$\,\mathcal{M}_{g,1}\,$ and $\,\mathcal{M}_{g,0}\,$ with the minimal 
possible number of twist generators. In $\,$\cite{Gervais}, the author 
gave a presentation considering either all possible Dehn twists, or 
just Dehn twists along non-separating curves. These two presentations 
appear to be very symmetric, but infinite. The aim of this article is 
to give a finite presentation of $\,\mathcal{M}_{g,n}$.

\vskip3mm\noindent
{\bf Notation.} Composition of diffeomorphisms in 
$\,\mathcal{M}_{g,n}\,$ will be written from right to left. For 
two elements $x$, $y$ of a multiplicative group, we will denote 
indifferently by $x^{-1}$ or $\bar{x}$ the inverse of $x$ and by 
$\,y(x)\,$ the conjugate $\,y\,x\,\bar{y}\,$ of $x$ by $y$.

\eject
Next, considering the curves of figure 1, we denote by 
$\mathcal{G}_{g,n}\,$ and $\mathcal{H}_{g,n}\,$ (we may on occasion omit the 
subscript ``$g,n$'' if there is no ambiguity) the following sets of
curves in $\,\Sigma_{g,n}$:
$$\begin{array}{rcl}
\mathcal{G}_{g,n} & = & \{\beta,\beta_{1},\ldots,\beta_{g-1},\alpha_{1},\ldots,
\alpha_{2g+n-2},(\gamma_{i,j})_{1\leq i,j\leq 2g+n-2,i\not =j}\,\},\\ &&\\
\mathcal{H}_{g,n} & = & \{\alpha_{1},\beta,\alpha_{2},\beta_{1},\gamma_{2,4},\beta_{2},\ldots,
\gamma_{2g-4,2g-2},\beta_{g-1},\gamma_{1,2},\\ &&\hspace*{45mm}
\alpha_{2g},\ldots,\alpha_{2g+n-2},\delta_{1},\ldots,\delta_{n-1}\,\}
\end{array}$$
where $\,\delta_{i}=\gamma_{2g-2+i,2g-1+i}\,$ is the i$^{\hbox{\scriptsize th}}$ 
boundary component. Note that $\,\mathcal{H}_{g,n}\,$ is a subset of
$\,\mathcal{G}_{g,n}$.

\vskip3mm
Finally, a triple $\,(i,j,k)\!\in\!\{1,\ldots,2g+n-2\}^{3}\,$ will 
be said to be {\em good} when:
$$\begin{array}{rl}
\hbox{i)} & (i,j,k)\!\not\in\!\bigl\{(x,x,x)\,/\,x\!\in\{1,\ldots,2g+n-2\}
            \bigr\},\\
\hbox{ii)} & i\leq j\leq k\ \hbox{ or }\ j\leq k\leq i\ \hbox{ or }\ 
             k\leq i\leq j\,.
\end{array}$$
\begin{center}
\courbes

figure 1
\end{center}

\vskip5mm\noindent
\begin{remark}
For $\,n=0\,$ or $\,n=1$, Wajnryb's generators are 
the Dehn twists relative to the curves of $\mathcal{H}$.
\end{remark}

We will give a presentation of $\,\mathcal{M}_{g,n}\,$ taking as 
generators the twists along the curves in $\mathcal{G}$. The 
relations will be of the following types.

\vskip3mm\noindent
{\bf The braids:} If $\alpha$ and $\beta$ are two curves in 
$\,\Sigma_{g,n}\,$ which do not intersect $\,$(resp. intersect in a 
single point), then the associated Dehn twists satisfy the 
relation $\,\tau_{\alpha}\tau_{\beta}=\tau_{\beta}\tau_{\alpha}\,$ (resp. 
$\,\tau_{\alpha}\tau_{\beta}\tau_{\alpha}=\tau_{\beta}\tau_{\alpha}\tau_{\beta}$).

\vskip3mm\noindent
{\bf The stars:} Concider a subsurface of $\,\Sigma_{g,n}\,$ which is 
homeomorphic to $\,\Sigma_{1,3}$. Then, if $\,\alpha_{1},\ 
\alpha_{2},\ \alpha_{3},\ \beta,\ \gamma_{1},\ 
\gamma_{2},\ \gamma_{3}\,$ are the curves described in figure 2, 
one has in $\,\mathcal{M}_{g,n}\,$ the relation
$$(\tau_{\alpha_{1}}\tau_{\alpha_{2}}\tau_{\alpha_{3}}
\tau_{\beta})^{3}=\tau_{\gamma_{1}}\tau_{\gamma_{2}}\tau_{\gamma_{3}}\,.$$
Note that if $\gamma_{3}\,$ bounds a disc in $\,\Sigma_{g,n}$, then 
this relation becomes
$$(\tau_{\alpha_{1}}\tau_{\alpha_{2}}\tau_{\alpha_{2}}
\tau_{\beta})^{3}=\tau_{\gamma_{1}}\tau_{\gamma_{2}}\,.$$

\begin{center}
\etoile

figure 2
\end{center}

\vskip5mm\noindent
{\bf The handles:} Pasting a cylinder on two boundary components of 
$\,\Sigma_{g-1,n+2}$, the twists along these two boundary curves 
become equal in $\,\Sigma_{g,n}$.

\vskip3mm
\begin{theorem}\label{principaltheorem}
For all $\,(g,n)\!\in\!{\mathbf{N}}^{\ast}\!\times\!{\mathbf{N}}$, the
mapping class group $\,\mathcal{M}_{g,n}\,$ admits a presentation with
generators $b,\,b_{_{1}},\ldots,b_{_{g-1}},a_{_{1}},\ldots,a_{_{2g+n-2}},$
$(c_{_{i,j}})_{1\leq i,j\leq 2g+n-2,\,i\not=j}\,$ and relations
\begin{list}{}{\labelsep=4mm\labelwidth=13mm\leftmargin=24mm\parsep=5mm\topsep=5mm}
\item[(A)] {\rm ``handles'':} $c_{_{2i,2i+1}}=c_{_{2i-1,2i}}\,$ for 
             all $i$, $\,1\leq i\leq g-1$,
\item[(T)] {\rm ``braids'':} for all $\,x,y\,$ among the generators, $xy=yx$
             if the associated curves are disjoint and $xyx\!=\!yxy$ if
             the associated curves intersect transversaly in a single point,
\item[(E$_{i,j,k}$)] {\rm ``stars'':} $\ c_{_{i,j}}c_{_{j,k}}c_{_{k,i}}
             =(a_{_{i}}a_{_{j}}a_{_{k}}b)^{3}\,$ for all good triples
             $\,(i,j,k)\,$, where $\,c_{_{l,l}}\!=\!1$.
\end{list}
\end{theorem}

\begin{remark}
It is clear that the handle relations are unnecessary: one has just to 
remove $\,c_{_{2,3}},\ldots,c_{_{2g-2,2g-1}}\,$ from 
$\,\mathcal{G}_{g,n}\,$ to eliminate them. But it is convenient for 
symmetry and notation to keep these generators.
\end{remark}

Let $\,G_{g,n}\,$ denote the group with presentation given by 
theorem~\ref{principaltheorem}. Since the set of generators for 
$\,G_{g,n}\,$ that we consider here is parametrized by 
$\,\mathcal{G}_{g,n}$, we will consider $\,\mathcal{G}_{g,n}\,$ as a 
subset of $\,G_{g,n}$. Consequently, $\,\mathcal{H}_{g,n}\,$ will also
be considered as a subset of $\,G_{g,n}$.

\hfill\break\indent
The paper is organized as follows. In section~\ref{Generators}, we prove
that $\,G_{g,n}\,$ is generated by $\,\mathcal{H}_{g,n}$.
Section~\ref{negalun} is devoted to the proof of theorem~\ref{principaltheorem}
when $\,n=1$. Finally, we conclude the proof in section~\ref{finpreuve} 
by proving that $\,G_{g,n}\,$ is isomorphic to $\,\mathcal{M}_{g,n}$.


\section{Generators for $G_{g,n}$ \label{Generators}}

\noindent
In this section, we prove the following proposition.

\begin{proposition} \label{generator}
$G_{g,n}\,$ is generated by $\,\mathcal{H}_{g,n}$.
\end{proposition}

\noindent
We begin by proving some relations in $G_{g,n}$.

\begin{lemma}\label{etoile}
For $\,i,j,k\!\in\! \{1,\ldots, 2g+n-2\}$, if $\,X_{1}=a_{_{i}}a_{_{j}},\ X_{2}=
bX_{1}b\,$ and $\,X_{3}=a_{_{k}}X_{2}a_{_{k}}$, then:
\begin{list}{(\roman{liste})}{\usecounter{liste}\parsep=2mm}
\item $X_{p}X_{q}=X_{q}X_{p}\,$ for all $\,p,q\!\in\!\{1,2,3\}$.
\item $(a_{_{i}}a_{_{j}}a_{_{k}}b)^{3}=X_{1}X_{2}X_{3}$,
\item $(a_{_{i}}a_{_{i}}a_{_{j}}b)^{3}=X_{1}^{2}X_{2}^{2}=
      (a_{_{i}}a_{_{j}}b)^{4}=(a_{_{i}}b\,a_{_{j}})^{4}$,
\item $a_{_{i}},\,a_{_{j}},\,a_{_{k}}\,$ and $b$ commute with
$\,(a_{_{i}}a_{_{j}}a_{_{k}}b)^{3}$.
\end{list}
\end{lemma}

\begin{remark}\label{rem}
Combining the braid relations
and lemma~\ref{etoile}, we get\linebreak[4] $\,(E_{i,j,k})=(E_{j,k,i})=
(E_{k,i,j})\,$ and $\,(E_{i,i,j})=(E_{i,j,j})$.
\end{remark}

\proof {\it (i) } Using relations {\it $\,$(T)$\,$}, one has
$$\begin{array}{rcl}
a_{_{i}}\,X_{2} & = & a_{_{i}}\,b\,a_{_{i}}\,a_{_{j}}\,b \\
&=& b\,a_{_{i}}\,b\,a_{_{j}}\,b \\
&=& b\,a_{_{i}}\,a_{_{j}}\,b\,a_{_{j}} \\
&=& X_{2}\,a_{_{j}}\,,
\end{array}$$
and in the same way, $\,a_{_{j}}\,X_{2}=X_{2}\,a_{_{i}}$. Thus, we 
get $\,X_{1}\,X_{2}=X_{2}\,X_{1}\,$ and $\,X_{1}\,X_{3}=X_{3}\,X_{1}\,$ since 
$\,X_{1}\,a_{_{k}}=a_{_{k}}\,X_{1}$.

\noindent
On the other hand, the braid relations imply
$$\begin{array}{rcl}
b(X_{3}) & = & 
b\,a_{_{k}}\,b\,a_{_{i}}\,a_{_{j}}\,b\,a_{_{k}}\,\bar{b} \\
&=& a_{_{k}}\,b\,a_{_{k}}\,a_{_{i}}\,a_{_{j}}\,\bar{a_{_{k}}}\,b\,a_{_{k}} \\
&=& X_{3}\,,
\end{array}$$
and we get $\,X_{2}\,X_{3}=X_{3}\,X_{2}$.

\vskip3mm\noindent
{\it (ii) } Using relations {\it $\,$(T)$\,$} and {\it $\,$(i)$\,$}, one 
obtains:
$$
\begin{array}{rcl}
X_{1}X_{2}X_{3} & = & X_{1}X_{3}X_{2} \\
&=& 
a_{_{i}}\,a_{_{j}}\,a_{_{k}}\,b\,a_{_{i}}\,a_{_{j}}\,b\,a_{_{k}}\,b\,
                                                      a_{_{i}}\,a_{_{j}}\,b\\
&=& 
a_{_{i}}\,a_{_{j}}\,a_{_{k}}\,b\,a_{_{i}}\,a_{_{j}}\,a_{_{k}}\,b\,a_{_{k}}\,
                                                       a_{_{i}}\,a_{_{j}}\,b\\
&=& (a_{_{i}}a_{_{j}}a_{_{k}}b)^{3}.
\end{array}
$$

\vskip3mm\noindent
{\it (iii) } Replacing $a_{_{k}}$ by $a_{_{i}}$ in $X_{3}$, we get
$$X_{3}=a_{_{i}}\,X_{2}\,a_{_{i}}=a_{_{i}}\,a_{_{j}}\,X_{2}=X_{1}\,X_{2}.$$
Thus, using relations {\it $\,$(T)}, {\it $\,$(i)$\,$} and {\it 
$\,$(ii)}, one has:
$$\begin{array}{rcl}
(a_{_{i}}a_{_{i}}a_{_{j}}b)^{3} & = & X_{1}X_{2}X_{1}X_{2}=X_{1}^{2}X_{2}^{2} \\
&=& a_{_{i}}\,a_{_{j}}\,b\,a_{_{i}}\,a_{_{j}}\,b\,a_{_{i}}\,
a_{_{j}}\,b\,a_{_{i}}\,a_{_{j}}\,b\,=\,(a_{_{i}}a_{_{j}}b)^{4} \\
&=& a_{_{i}}\,b\,a_{_{j}}\,b\,a_{_{i}}\,b\,
a_{_{j}}\,b\,a_{_{i}}\,b\,a_{_{j}}\,b  \\
&=& a_{_{i}}\,b\,a_{_{j}}\,a_{_{i}}\,b\,a_{_{i}}\,
a_{_{j}}\,b\,a_{_{i}}\,a_{_{j}}\,b\,a_{_{j}}  \\
&=&(a_{_{i}}b\,a_{_{j}})^{4}.
\end{array}$$

 %
 %
 
\vskip3mm\noindent
{\it	(iv) } One has just to apply the star and braid relations.
\eproof

\begin{lemma} \label{lantern}
For all good triples $\,(i,j,k)$, one has in $\,G_{g,n}\,$ the relation
$$(L_{_{i,j,k}})\ \ 
a_{_{i}}\,c_{_{i,j}}\,c_{_{j,k}}\,a_{_{k}}=c_{_{i,k}}\,a_{_{j}}\,X\,a_{_{j}}\,
\bar{X}=c_{_{i,k}}\,\bar{X}\,a_{_{j}}\,X\,a_{_{j}}$$
where $\,X\!=\!b\,a_{_{i}}\,a_{_{k}}\,b$.
\end{lemma}

\begin{remark}
These relations are just the well known {\em lantern} relations.
\end{remark}

\proof If $\,X_{1}\!=\!a_{_{i}}\,a_{_{k}}\,$ and 
$\,X_{3}\!=\!a_{_{j}}\,X\,a_{_{j}}\,$, one has by lemma~\ref{etoile} and the
star relations $\,(E_{_{i,j,k}})\,$ and $\,(E_{_{i,k,k}})\,$:
$$X_{1}\,X\,X_{3}=c_{_{i,j}}\,c_{_{j,k}}\,c_{_{k,i}}\ \hbox{ and }\ X_{1}^{2}\,
X^{2}=c_{_{i,k}}\,c_{_{k,i}}\,.$$

\noindent
From this, we get, using the braid relations, that
$$\bar{c_{_{k,i}}}\,X_{1}\,X=c_{_{i,j}}\,c_{_{j,k}}\,\bar{X_{3}}=c_{_{i,k}}\,
\bar{X}\,\bar{X_{1}}\,,$$
that is to say, by lemma~\ref{etoile} and {\it (T)},
$$a_{_{i}}\,c_{_{i,j}}\,c_{_{j,k}}\,a_{_{k}}=c_{_{i,k}}\,\bar{X}\,a_{_{j}}\,X\,
a_{_{j}}=c_{_{i,k}}\,a_{_{j}}\,X\,a_{_{j}}\,\bar{X}\,.$$
\eproof

\begin{lemma}\label{ak}
For all $i,k$ such that $\,1\!\leq i\leq g-1\,$ and 
$\,k\not=2i-1,2i$, one has in $\,G_{g,n}$
$$a_{_{k}}\,=\,b\,a_{_{2i}}\,b_{_{i}}\,a_{_{2i-1}}\,b\,
\bar{c_{_{2i,2i-1}}}\,a_{_{2i}}\,c_{_{2i,k}}(b_{_{i}})\,.$$ 
\end{lemma}

\vskip3mm
\proof If $\,X\!=\!b\,a_{_{2i-1}}\,a_{_{2i}}\,b$, one has by the 
lantern relations
$$(L_{2i,k,2i-1}):\,\ a_{_{2i}}\,c_{_{2i,k}}\,c_{_{k,2i-1}}\,a_{_{2i-1}}=
c_{_{2i,2i-1}}\,\bar{X}\,a_{_{k}}\,X\,a_{_{k}}\,,$$
which implies
$$\bar{c_{_{2i,2i-1}}}\,a_{_{2i}}\,c_{_{2i,k}}=
\bar{X}\,a_{_{k}}\,X\,a_{_{k}}\,\bar{a_{_{2i-1}}}\,\bar{c_{_{k,2i-1}}}\,.$$

\vskip3mm\noindent
Thus, denoting $\,b\,a_{_{2i}}\,b_{_{i}}\,a_{_{2i-1}}\,b\,\bar{c_{_{2i,2i-1}}}\,
a_{_{2i}}\,c_{_{2i,k}}(b_{_{i}})\,$ by $y$,  we can compute using the
relations {\it $\,$(T)}: 

$$\begin{array}{rcl}
y & = & b\,a_{_{2i}}\,b_{_{i}}\,a_{_{2i-1}}\,b\,\bar{X}\,a_{_{k}}\,X\,a_{_{k}}\,
        \bar{a_{_{2i-1}}}\,\bar{c_{_{k,2i-1}}}(b_{_{i}}) \\
&=& b\,a_{_{2i}}\,b_{_{i}}\,a_{_{2i-1}}\,b\,\bar{b}\,\bar{a_{_{2i-1}}}\,
    \bar{a_{_{2i}}}\,\bar{b}\,a_{_{k}}\,b\,a_{_{2i-1}}\,a_{_{2i}}\,b\,
    (b_{_{i}}) \\
&=& b\,\bar{b_{_{i}}}\,a_{_{2i}}\,b_{_{i}}\,a_{_{k}}\,b\,\bar{a_{_{k}}}\,
    \bar{b_{_{i}}}\,(a_{_{2i}}) \\
&=& b\,a_{_{k}}\,\bar{b_{_{i}}}\,a_{_{2i}}\,\bar{a_{_{2i}}}(b) \\
&=& b\,\bar{b}(a_{_{k}}) \\
&=& a_{_{k}}.
\end{array}
$$
\eproof
\vskip3mm\noindent
{\bf Proof of proposition~\ref{generator}.\ \,}If $H$ denotes the subgroup
of $\,G_{g,n}\,$ generated by $\,\mathcal{H}_{g,n}\,$, we have to 
prove that $\,\mathcal{G}_{g,n}\!\subset\! H$.

\vskip3mm\noindent
a) We first prove inductively that $\,a_{_{2i-1}},\,a_{_{2i}},\,c_{_{2i-1,2i}}\,$
and $\,c_{_{2i,2i-1}}\,$ are elements of $H$ for all $i$, $\,1\leq i\leq 
g-1$.

For $\,i\!=\!1$, one obtains $\,a_{_{1}},\,a_{_{2}}\,$ and 
$\,c_{_{1,2}}\,$ which are in $H$, and the 
relation $\,(E_{1,2,2})\,$ gives $\,c_{_{2,1}}\!=\!(a_{_{1}}a_{_{2}}
a_{_{2}}b)^{3}\bar{c_{_{1,2}}}\in H$. So, suppose inductively that
$\,a_{_{2i-1}},\,a_{_{2i}},\,c_{_{2i-1,2i}},\,c_{_{2i,2i-1}}\,$ are elements
of $H\,$ ($i\leq g-2$) and let us prove that $\,a_{_{2i+1}},\,a_{_{2i+2}},\,
c_{_{2i+1,2i+2}},\,c_{_{2i+2,2i+1}}\,$ are also in $H$. Recall that by 
the handle relations, one has $\,c_{_{2i,2i+1}}\!=\!
c_{_{2i-1,2i}}\!\in\! H$. Applying lemma~\ref{ak} respectively with 
$\,k\!=\!2i+1\,$ and $\,k\!=\!2i+2$, we obtain

$$a_{_{2i+1}}\,=\,b\,a_{_{2i}}\,b_{_{i}}\,a_{_{2i-1}}\,b\,
\bar{c_{_{2i,2i-1}}}\,a_{_{2i}}\,c_{_{2i,2i+1}}(b_{_{i}})\in H\,,$$
$$a_{_{2i+2}}\,=\,b\,a_{_{2i}}\,b_{_{i}}\,a_{_{2i-1}}\,b\,
\bar{c_{_{2i,2i-1}}}\,a_{_{2i}}\,c_{_{2i,2i+2}}(b_{_{i}})\in H\,.$$

\vskip4mm\noindent
The star relations allow us to conclude the induction as follows:

$$(E_{_{2i,2i+2,2i+2}})\,:\ \ \  
c_{_{2i,2i+2}}\,c_{_{2i+2,2i}}=(a_{_{2i}}\,a_{_{2i+2}}\,b)^{4},$$
which gives $\,c_{_{2i+2,2i}}\!\in\! H\,$ 
($\gamma_{2i,2i+2}\!\in\!\mathcal{H}_{g,n}\,$ by definition);

$$(E_{_{2i,2i+1,2i+2}})\,:\ \ \  
c_{_{2i,2i+1}}c_{_{2i+1,2i+2}}c_{_{2i+2,2i}}=
(a_{_{2i}}a_{_{2i+1}}a_{_{2i+2}}b)^{3},$$
which gives $\,c_{_{2i+1,2i+2}}\!\!\in\! H$;

$$(E_{_{2i+1,2i+2,2i+2}})\,:\ \ \  
c_{_{2i+1,2i+2}}\,c_{_{2i+2,2i+1}}=
(a_{_{2i+1}}\,a_{_{2i+2}}\,b)^{4},$$
which gives $\,c_{_{2i+2,2i+1}}\!\in\! H$.

\vskip7mm\noindent
b) By lemma~\ref{ak}, one has ($i=g-1$ and $k=2g-1$)
$$a_{_{2g-1}}=b\,a_{_{2g-2}}\,b_{_{g-1}}\,a_{_{2g-3}}\,b\,
\bar{c_{_{2g-2,2g-3}}}\,a_{_{2g-2}}\,c_{_{2g-2,2g-1}}(b_{_{g-1}}).$$
Recall that $\,c_{_{2g-2,2g-1}}\!=\!c_{_{2g-3,2g-2}}\!\in\! H$. Thus,
combined with the case a), this relation implies $\,a_{_{2g-1}}\!\in\! H$.

\vskip7mm\noindent
c) It remains to prove that $\,c_{_{i,j}}\in H\,$ for all $i,j$.

\vskip3mm
$\ast\ $ By definition of $H$ and the case a), one has $\,c_{_{i,i+1}}\in H\,$
for all $i$ such that $\,1\leq i\leq 2g+n-3$.

\vskip3mm
$\ast\ $ Let us show that $\,c_{_{1,j}}\,$ and $\,c_{_{j,1}}\,$ are 
elements of $H$ for all $j$ such that $\,2\leq j\leq 2g+n-2$.

We have already seen that $\,c_{_{1,2}},\,c_{_{2,1}}\!\in\! H$.
Thus, suppose inductively that $\,c_{_{1,j}},c_{_{j,1}}\in H\,$ ($j\leq 2g+n-3$).
Using the star relations, one obtains:

\vskip3mm\noindent
\begin{center}
$(E_{_{1,j,j+1}})\,$: $\,c_{_{1,j}}\,c_{_{j,j+1}}\,c_{_{j+1,1}}=
(a_{_{1}}\,a_{_{j}}\,a_{_{j+1}}\,b)^{3}$, which gives $\,c_{_{j+1,1}}\in H$,

\vskip3mm\noindent
$(E_{_{1,j+1,j+1}})\,$: $\,c_{_{1,j+1}}\,c_{_{j+1,1}}=
(a_{_{1}}\,a_{_{j+1}}\,b)^{4}$, which gives $\,c_{_{1,j+1}}\in H$.
\end{center}

\vskip3mm
$\ast\ $ Now, fix $j$ such that $\,2\!\leq\! j\!\leq\! 2g+n-2\,$ and let 
us show that $\,c_{_{i,j}},c_{_{j,i}}\in H\,$ for all $i$, $\,1\leq 
i<j$. Once more, the star relations allow us to prove this using an 
inductive argument:

\vskip3mm\noindent
\begin{center}
$(E_{_{i,i+1,j}})\,$: $\,c_{_{i,i+1}}\,c_{_{i+1,j}}\,c_{_{j,i}}=
(a_{_{i}}\,a_{_{i+1}}\,a_{_{j}}\,b)^{3}$, which gives $\,c_{_{i+1,j}}\in H$,

\vskip3mm\noindent
$(E_{_{i+1,j,j}})\,$: $\,c_{_{i+1,j}}\,c_{_{j,i+1}}=
(a_{_{i+1}}\,a_{_{j}}\,b)^{4}$, which gives $\,c_{_{j,i+1}}\in H$.
\end{center}

\eproof


\section{Proof of theorem~\ref{principaltheorem} for $n=1$ \label{negalun}}

\noindent
Let us recall Wajnryb's result:

\begin{theorem}[\cite{Wajnryb}]
$\ \mathcal{M}_{g,1}\,$ admits a presentation with generators\linebreak[4]
$\,\{\tau_{\alpha}\,/\alpha\!\in\!\mathcal{H}\}\,$
and  relations

\begin{list}{(\Roman{liste})}{\usecounter{liste}\labelwidth=8mm
                              \labelsep=4mm\leftmargin=17mm\itemsep3mm}
\item $\,\tau_{\lambda}\tau_{\mu}\tau_{\lambda}=\tau_{\mu}\tau_{\lambda}\tau_{\mu}\,$
      if $\lambda$ and $\mu$ intersect transversaly in a single point, and
      $\,\tau_{\lambda}\tau_{\mu}=\tau_{\mu}\tau_{\lambda}\,$ if $\lambda$ and $\mu$
      are disjoint.
\item $(\tau_{\alpha_{1}}\tau_{\beta}\tau_{\alpha_{2}})^{4}=\tau_{\gamma_{1,2}}\,\theta\,$
      where $\,\theta=\tau_{\beta_{1}}\tau_{\alpha_{2}}\tau_{\beta}\tau_{\alpha_{1}}
      \tau_{\alpha_{1}}\tau_{\beta}\tau_{\alpha_{2}}\tau_{\beta_{1}}(\tau_{\gamma_{1,2}})$.
\item $\tau_{\alpha_{2}}\tau_{\alpha_{1}}\varphi\,\tau_{\gamma_{2,4}}=
       \bar{t_{1}}\,\bar{t_{2}}\,\tau_{\gamma_{1,2}}\,t_{2}\,t_{1}\,\bar{t_{2}}\,
      \tau_{\gamma_{1,2}}\,t_{2}\,\tau_{\gamma_{1,2}}\,\ $ where 
      \begin{center}
      $\,t_{1}=\tau_{\beta}\tau_{\alpha_{1}}
      \tau_{\alpha_{2}}\tau_{\beta}\,$, $\,t_{2}=\tau_{\beta_{1}}\tau_{\alpha_{2}}
      \tau_{\gamma_{2,4}}\tau_{\beta_{1}}\,$,
      
      \noindent $\,\varphi=\tau_{\beta_{2}}\tau_{\gamma_{2,4}}
      \tau_{\beta_{1}}\tau_{\alpha_{2}}\tau_{\beta}\,\sigma(\omega)$,
      $\,\sigma=\bar{\tau_{\gamma_{2,4}}}\,\bar{\tau_{\beta_{2}}}\,\bar{t_{2}}
      (\tau_{\gamma_{1,2}})\,$
      
      \noindent
      and  $\,\omega=\bar{\tau_{\alpha_{1}}}\,\bar{\tau_{\beta}}\,
      \bar{\tau_{\alpha_{2}}}\,\bar{\tau_{\beta_{1}}}(\tau_{\gamma_{1,2}})$.
      \end{center}
\end{list}
\end{theorem}

\begin{remark}
When $\,g\!=\!1$, one just needs the relations $\,${\it (I)}. The relations 
$\,${\it (II)}$\,$ and $\,${\it (III)}$\,$ appear respectively for $\,g\!=\!2\,$ 
and $\,g\!=\!3$.
\end{remark}

\vskip3mm
Denote by $\,\Phi\!:\!G_{g,1}\!\rightarrow\! \mathcal{M}_{g,1}\,$ the map
which associates to each generator $a$ of $\,G_{g,1}\,$ the corresponding 
twist $\tau_{\alpha}$. Since the relations {\it (A), (T)} and {\it 
(E$_{i,j,k}$)} are satisfied in $\,\mathcal{M}_{g,1}$, $\Phi$ is an homomorphism.

\noindent
Now, consider $\,\Psi:\mathcal{M}_{g,1}\rightarrow G_{g,1}\,$ defined by 
$\,\Psi(\tau_{\alpha})=a\,$ for all 
$\,\alpha\in\mathcal{H}$.

\begin{lemma}\label{psi}
$\Psi$ is an homomorphism.
\end{lemma}

This lemma allows us to prove the theorem~\ref{principaltheorem} for 
$\,n=1$. Indeed, since $\,\mathcal{M}_{g,1}$ is generated by
$\,\{\tau_{\alpha}\,/\,\alpha\in\mathcal{H}_{g,1}\}$, 
one has $\,\Phi\circ\Psi=Id_{_{\mathcal{M}_{g,1}}}$. On the other hand, 
$\,\{a\,/\,\alpha\in\mathcal{H}_{g,1}\}\,$ 
generates $\,G_{g,1}\,$ by proposition~\ref{generator}, so
$\,\Psi\circ\Phi=Id_{_{G_{g,1}}}$.

\vskip5mm\noindent
{\bf Proof of lemma~\ref{psi}.\,} We have to show that the 
relations {\it (I)}, {\it (II)} and {\it (III)} are satisfied in 
$G_{g,1}$. Relations {\it (I)} are braid relations and are therefore 
satisfied by {\it (T)}. Let us look at the relation {\it (II)}. The 
star relation $\,(E_{_{1,2,2}})$, together with
lemma~\ref{etoile}, gives $\,(a_{_{1}}\,b\,a_{_{2}})^4=c_{_{1,2}}\,c_{_{2,1}}$.
Thus, relation {\it (II)} is satisfied in $G_{g,1}\,$ if and only 
if $\,\Psi(\theta)=c_{_{2,1}}$. Let us compute:
$$\begin{array}{rcll}
\Psi(\theta) & = & 
b_{_{1}}\,a_{_{2}}\,b\,a_{_{1}}\,a_{_{1}}\,b\,a_{_{2}}\,b_{_{1}}(c_{_{1,2}}) & \\
&=&b_{_{1}}\,a_{_{2}}\,b\,a_{_{1}}\,a_{_{1}}\,b\,a_{_{2}}\,\bar{c_{_{1,2}}}
(b_{_{1}}) & \hbox{by {\it (T)}}, \\
&=&b_{_{1}}\,a_{_{2}}\,b\,a_{_{1}}\,a_{_{1}}\,b\,a_{_{2}}\,(\bar{a_{_{1}}}\,
\bar{a_{_{1}}}\,\bar{a_{_{2}}}\,\bar{b})^{3}c_{_{2,1}}(b_{_{1}}) & 
\hbox{by } \,(E_{_{1,1,2}}), \\
&=&b_{_{1}}\,\bar{b}\,\bar{a_{_{1}}}\,\bar{a_{_{1}}}\,\bar{b}\,\bar{a_{_{1}}}\,
\bar{a_{_{1}}}\,c_{_{2,1}}(b_{_{1}}) & \hbox{by lemma~\ref{etoile}}, \\
&=&b_{_{1}}\,\bar{b_{_{1}}}(c_{_{2,1}}) & \hbox{by {\it (T)}}, \\
&=& c_{_{2,1}}. &
\end{array}$$

\vskip3mm\noindent
Wajnryb's relation {\it (III)} is nothing but a lantern relation. 
Via $\Psi$, it becomes in $\,G_{g,1}\,$
$$a_{_{2}}\,a_{_{1}}\,f\,c_{_{2,4}}=l\,m\,c_{_{1,2}}\ \ \ (\ast)$$

\vskip3mm\noindent
where 
$\,m=\bar{b_{_{1}}}\,\bar{a_{_{2}}}\,\bar{c_{_{2,4}}}\,\bar{b_{_{1}}}
                                                               (c_{_{1,2}})$,
$\,l=\bar{b}\,\bar{a_{_{1}}}\,\bar{a_{_{2}}}\,\bar{b}(m)\,$ and 
$\,f=b_{_{2}}\,c_{_{2,4}}\,b_{_{1}}\,a_{_{2}}\,b\,s(w)$, with 
$\,s=\Psi(\sigma)=\bar{c_{_{2,4}}}\,\bar{b_{_{2}}}(m)\,$ and 
$\,w=\Psi(\omega)=\bar{a_{_{1}}}\,\bar{b}\,\bar{a_{_{2}}}\,\bar{b_{_{1}}}
                                                               (c_{_{1,2}})$.

\vskip3mm\noindent
In $G_{g,1}$, the lantern relation $\,(L_{_{1,2,4}})\,$ yields
$$a_{_{1}}\,c_{_{1,2}}\,c_{_{2,4}}\,a_{_{4}}=
c_{_{1,4}}\,\bar{X}\,a_{_{2}}\,X\,a_{_{2}}\ \ \ (L_{_{1,2,4}})$$
where $\,X=b\,a_{_{1}}\,a_{_{4}}\,b$. To prove that the relation $\,(\ast)\,$ 
is satisfied in $\,G_{g,1}$, we will see that it is exactly the 
conjugate of the relation $\,(L_{_{1,2,4}})\,$ by 
$\,h=b_{_{2}}\,a_{_{4}}\,\bar{c_{_{4,1}}}\,\bar{b_{_{2}}}\,b
\,a_{_{2}}\,a_{_{1}}\,b\,b_{_{1}}\,c_{_{1,2}}\,a_{_{2}}\,b_{_{1}}$. 
This will be done by proving the following seven equalities in 
$\,G_{g,1}\,$:

$$\begin{array}{c}
\hbox{ 1) }\,h(a_{_{1}})=a_{_{2}}\ \ \ \hbox{ 2) }\,h(c_{_{1,2}})=a_{_{1}}
\ \ \ \hbox{ 3) }\,h(c_{_{2,4}})=f\ \ \ \hbox{ 4) }\,h(a_{_{4}})=c_{_{2,4}} \\ \\
\hbox{ 5) }\,h(c_{_{1,4}})=l\ \ \ \hbox{ 6) }\,h(a_{_{2}})=c_{_{1,2}}
\ \ \ \hbox{ 7) }\,h\bar{X}(a_{_{2}})=m.
\end{array}$$

\vskip3mm\noindent
1) Just applying the relations {\it (T)}, one obtains:
$$\begin{array}{rcl}
h(a_{_{1}}) & = & b_{_{2}}\,a_{_{4}}\,\bar{c_{_{4,1}}}\,\bar{b_{_{2}}}\,b
\,a_{_{2}}\,a_{_{1}}\,b\,b_{_{1}}\,c_{_{1,2}}\,a_{_{2}}\,b_{_{1}}(a_{_{1}}) \\
&=& b_{_{2}}\,a_{_{4}}\,\bar{c_{_{4,1}}}\,\bar{b_{_{2}}}\,b
\,a_{_{2}}\,a_{_{1}}\,\bar{a_{_{1}}}(b) \\
&=& b_{_{2}}\,a_{_{4}}\,\bar{c_{_{4,1}}}\,
\bar{b_{_{2}}}\,b\,\bar{b}(a_{_{2}}) \\
&=& a_{_{2}}\,.
\end{array}$$

\noindent
2) Using the relations {\it (T)} again, we get
$$
\begin{array}{rcl}
h(c_{_{1,2}}) & = & b_{_{2}}\,a_{_{4}}\,\bar{c_{_{4,1}}}\,\bar{b_{_{2}}}\,b
\,a_{_{2}}\,a_{_{1}}\,b\,b_{_{1}}\,c_{_{1,2}}\,a_{_{2}}\,b_{_{1}}(c_{_{1,2}}) \\
&=& b_{_{2}}\,a_{_{4}}\,\bar{c_{_{4,1}}}\,\bar{b_{_{2}}}\,b\,a_{_{2}}\,
a_{_{1}}\,b\,b_{_{1}}\,c_{_{1,2}}\,a_{_{2}}\,\bar{c_{_{1,2}}}(b_{_{1}})  \\
&=& b_{_{2}}\,a_{_{4}}\,\bar{c_{_{4,1}}}\,\bar{b_{_{2}}}\,b
\,a_{_{2}}\,a_{_{1}}\,b\,b_{_{1}}\,\bar{b_{_{1}}}(a_{_{2}})  \\
&=& b_{_{2}}\,a_{_{4}}\,\bar{c_{_{4,1}}}\,\bar{b_{_{2}}}\,b
\,a_{_{2}}\,a_{_{1}}\,\bar{a_{_{2}}}(b)  \\
&=& b_{_{2}}\,a_{_{4}}\,\bar{c_{_{4,1}}}\,\bar{b_{_{2}}}\,b
\,\bar{b}(a_{_{1}})  \\
&=& a_{_{1}}\,.
\end{array}
$$

\noindent
3) The relation $\,(L_{_{2,3,4}})\,$ yields
$$a_{_{2}}\,c_{_{2,3}}\,c_{_{3,4}}\,a_{_{4}}=
  c_{_{2,4}}\,\bar{Y}\,a_{_{3}}\,Y\,a_{_{3}}\ \ \ \hbox{ where }\ 
  Y=b\,a_{_{2}}\,a_{_{4}}\,b.$$
Since $\,c_{_{2,3}}\!\!=\!\!c_{_{1,2}}\,$ by the handle relations, this 
equality implies the following one:
$$\bar{c_{_{2,4}}}\,a_{_{2}}\,c_{_{1,2}}= \bar{Y}\,a_{_{3}}\,Y\,a_{_{3}}\,
\bar{a_{_{4}}}\,\bar{c_{_{3,4}}}\ \ \ \ \ (1).$$
  
\noindent
From this, we get:
$$\begin{array}{rcll}
h(c_{_{2,4}}) & = & b_{_{2}}\,a_{_{4}}\,\bar{c_{_{4,1}}}\,\bar{b_{_{2}}}\,b\,
a_{_{2}}\,a_{_{1}}\,b\,b_{_{1}}\,c_{_{1,2}}\,a_{_{2}}\,b_{_{1}}(c_{_{2,4}}) & \\

&=& b_{_{2}}\,a_{_{4}}\,\bar{c_{_{4,1}}}\,\bar{b_{_{2}}}\,b\,a_{_{2}}\,
a_{_{1}}\,b\,b_{_{1}}\,\bar{c_{_{2,4}}}\,c_{_{1,2}}\,a_{_{2}}(b_{_{1}})
& \hbox{by {\it (T)}} \\

&=& b_{_{2}}\,a_{_{4}}\,\bar{c_{_{4,1}}}\,\bar{b_{_{2}}}\,b\,a_{_{2}}\,
a_{_{1}}\,b\,b_{_{1}}\,\bar{Y}\,a_{_{3}}\,Y\,a_{_{3}}\,\bar{a_{_{4}}}\,
\bar{c_{_{3,4}}}(b_{_{1}}) & \hbox{by }\,(1) \\

&=& b_{_{2}}\,a_{_{4}}\,\bar{c_{_{4,1}}}\,\bar{b_{_{2}}}\,b\,a_{_{2}}\,
a_{_{1}}\,b\,b_{_{1}}\,\bar{b}\,\bar{a_{_{2}}}\,\bar{a_{_{4}}}\,\bar{b}\,
a_{_{3}}\,b\,a_{_{2}}\,a_{_{4}}\,b(b_{_{1}}) & \hbox{by {\it (T)}} \\

&=& b_{_{2}}\,a_{_{4}}\,\bar{c_{_{4,1}}}\,\bar{b_{_{2}}}\,b\,a_{_{1}}\,
\bar{b_{_{1}}}\,a_{_{2}}\,b_{_{1}}\,\bar{a_{_{4}}}\,a_{_{3}}\,b\,
\bar{a_{_{3}}}\,\bar{b_{_{1}}}(a_{_{2}}) & \hbox{by {\it (T)}} \\

&=& b_{_{2}}\,a_{_{4}}\,\bar{c_{_{4,1}}}\,\bar{b_{_{2}}}\,b\,a_{_{1}}\,
\bar{b_{_{1}}}\,a_{_{2}}\,\bar{a_{_{4}}}\,a_{_{3}}\,\bar{a_{_{2}}}(b)
 & \hbox{by {\it (T)}} \\

&=& b_{_{2}}\,a_{_{4}}\,\bar{c_{_{4,1}}}\,b\,a_{_{1}}\,
a_{_{3}}\,\bar{b_{_{2}}}\,b(a_{_{4}}) & \hbox{by {\it (T)}} \\

&=& b_{_{2}}\,a_{_{4}}\,(\bar{a_{_{1}}}\,\bar{a_{_{3}}}\,\bar{a_{_{4}}}\,
\bar{b})^{3}\,c_{_{1,3}}\,c_{_{3,4}}\,b\,a_{_{1}}\,a_{_{3}}\,
\bar{b_{_{2}}}\,b(a_{_{4}}) & \hbox{by }\,(E_{_{1,3,4}}) \\

&=& b_{_{2}}\,\bar{a_{_{1}}}\,\bar{a_{_{3}}}\,\bar{b}\,
(\bar{a_{_{1}}}\,\bar{a_{_{3}}}\,\bar{a_{_{4}}}\,\bar{b})^{2}\,
b\,a_{_{1}}\,a_{_{3}}\,c_{_{3,4}}\,b\,a_{_{4}}(b_{_{2}}) 
& \hbox{by {\it (T)}} \\

&=& b_{_{2}}\,\bar{a_{_{1}}}\,\bar{a_{_{3}}}\,\bar{b}\,
\bar{a_{_{1}}}\,\bar{a_{_{3}}}\,\bar{b}\,\bar{a_{_{4}}}\,\bar{b}\,
b\,a_{_{4}}\,\bar{b_{_{2}}}(c_{_{3,4}}) 
& \hbox{by {\it (T)}} \\

&=& c_{_{3,4}} & \hbox{by {\it (T)}.} \\
\end{array}$$
Now, if $\,x\!=\!c_{_{1,2}}\,b_{_{1}}\,c_{_{2,4}}\,a_{_{2}}\,b_{_{1}}\,
b_{_{2}}\,c_{_{2,4}}\,\bar{a_{_{1}}}\,\bar{b}\,\bar{a_{_{2}}}\,
\bar{b_{_{1}}}(c_{_{1,2}})$, one has
$$f=b_{_{2}}\,c_{_{2,4}}\,b_{_{1}}\,a_{_{2}}\,b\,\bar{c_{_{2,4}}}\,
\bar{b_{_{2}}}\,\bar{b_{_{1}}}\,\bar{a_{_{2}}}\,\bar{c_{_{2,4}}}\,
\bar{b_{_{1}}}(x)\,.$$
First, let us compute $x$:
$$\begin{array}{rcll}
x & = & c_{_{1,2}}\,b_{_{1}}\,c_{_{2,4}}\,a_{_{2}}\,b_{_{1}}\,
b_{_{2}}\,c_{_{2,4}}\,\bar{a_{_{1}}}\,\bar{b}\,\bar{a_{_{2}}}\,
\bar{b_{_{1}}}(c_{_{1,2}}) &  \\

&=& c_{_{1,2}}\,b_{_{1}}\,c_{_{2,4}}\, a_{_{2}}\,b_{_{1}}\,
b_{_{2}}\,c_{_{2,4}}\,c_{_{1,2}}\,\bar{a_{_{1}}}\,\bar{b}\,
\bar{a_{_{2}}}(b_{_{1}}) & \hbox{by {\it (T)}} \\

&=& c_{_{1,2}}\,b_{_{1}}\,c_{_{2,4}}\, a_{_{2}}\,b_{_{1}}\,
b_{_{2}}\,(a_{_{1}}\,a_{_{2}}\,a_{_{4}}\,b)^{3}\,\bar{c_{_{4,1}}}\,
\bar{a_{_{1}}}\,\bar{b}\,
\bar{a_{_{2}}}(b_{_{1}}) & \hbox{by }\,(E_{_{1,2,4}}) \\

&=& c_{_{1,2}}\,b_{_{1}}\,c_{_{2,4}}\, a_{_{2}}\,b_{_{1}}\,
b_{_{2}}\,(a_{_{1}}\,a_{_{2}}\,a_{_{4}}\,b)^{2}\,a_{_{1}}\,a_{_{2}}\,
a_{_{4}}\,b\,\bar{a_{_{1}}}\,\bar{b}\,
\bar{a_{_{2}}}(b_{_{1}}) & \hbox{by {\it (T)}} \\

&=& c_{_{1,2}}\,b_{_{1}}\,c_{_{2,4}}\, a_{_{2}}\,b_{_{1}}\,
b_{_{2}}\,(a_{_{1}}\,a_{_{2}}\,a_{_{4}}\,b)^{2}\,a_{_{4}}\,a_{_{2}}\,
\bar{b}\,a_{_{1}}\,b\,\bar{b}\,
\bar{a_{_{2}}}(b_{_{1}}) & \hbox{by {\it (T)}} \\

&=& c_{_{1,2}}\,b_{_{1}}\,c_{_{2,4}}\, a_{_{2}}\,b_{_{1}}\,
b_{_{2}}\,(a_{_{1}}\,a_{_{2}}\,a_{_{4}}\,b)^{2}\,a_{_{4}}\,
\bar{b}\,\bar{a_{_{2}}}\,b(b_{_{1}}) & \hbox{by {\it (T)}} \\

&=& c_{_{1,2}}\,b_{_{1}}\,c_{_{2,4}}\, a_{_{2}}\,b_{_{1}}\,
b_{_{2}}\,a_{_{1}}\,a_{_{2}}\,a_{_{4}}\,b\,a_{_{1}}\,a_{_{2}}\,b\,
a_{_{4}}\,b\,\bar{b}\,\bar{a_{_{2}}}(b_{_{1}}) & \hbox{by {\it (T)}} \\

&=& c_{_{1,2}}\,b_{_{1}}\,c_{_{2,4}}\, a_{_{2}}\,b_{_{1}}\,
b_{_{2}}\,a_{_{1}}\,a_{_{2}}\,a_{_{4}}\,b\,a_{_{1}}\,\bar{b}\,a_{_{2}}\,
b(b_{_{1}}) & \hbox{by {\it (T)}} \\

&=& c_{_{1,2}}\,b_{_{1}}\,c_{_{2,4}}\, b_{_{2}}\,a_{_{2}}\,b_{_{1}}\,
a_{_{2}}\,a_{_{4}}\,b\,a_{_{1}}\,b\,\bar{b}\,a_{_{2}}
(b_{_{1}}) & \hbox{by {\it (T)}} \\

&=& c_{_{1,2}}\,b_{_{1}}\,c_{_{2,4}}\, b_{_{2}}\,b_{_{1}}\,a_{_{2}}\,b_{_{1}}\,
a_{_{4}}\,b\,\bar{b_{_{1}}}(a_{_{2}}) & \hbox{by {\it (T)}} \\

&=& c_{_{1,2}}\,b_{_{1}}\,c_{_{2,4}}\, b_{_{2}}\,b_{_{1}}\,a_{_{2}}\,
a_{_{4}}\,\bar{a_{_{2}}}(b) & \hbox{by {\it (T)}} \\

&=& c_{_{1,2}}\,b_{_{1}}\,c_{_{2,4}}\, b_{_{2}}\,\bar{b}(a_{_{4}})
 & \hbox{by {\it (T)}.} \\
\end{array}$$
Next, using the braid relations, we prove that $\,b_{_{1}},\ c_{_{2,4}},
\ b_{_{2}}\,$ and $\,a_{_{2}}\,$ commute with $x$ :
$$b_{_{1}}(x) = b_{_{1}}\,c_{_{1,2}}\,b_{_{1}}\,c_{_{2,4}}\, b_{_{2}}\,
\bar{b}(a_{_{4}}) = c_{_{1,2}}\,b_{_{1}}\,c_{_{1,2}}\,c_{_{2,4}}\, b_{_{2}}\,
    \bar{b}(a_{_{4}}) = x,$$
$$c_{_{2,4}}(x) = c_{_{1,2}}\,b_{_{1}}\,c_{_{2,4}}\,b_{_{1}}\, b_{_{2}}\,
\bar{b}(a_{_{4}}) = x,$$
$$b_{_{2}}(x) = c_{_{1,2}}\,b_{_{1}}\,b_{_{2}}\,c_{_{2,4}}\,b_{_{2}}\,
\bar{b}(a_{_{4}}) = c_{_{1,2}}\,b_{_{1}}\,c_{_{2,4}}\,b_{_{2}}\,c_{_{2,4}}\,
\bar{b}(a_{_{4}}) = x,$$
$$\begin{array}{rcll}
a_{_{2}}(x) & = & a_{_{2}}\,c_{_{1,2}}\,b_{_{1}}\,c_{_{2,4}}\,a_{_{2}}\,
b_{_{1}}\,b_{_{2}}\,c_{_{2,4}}\,\bar{a_{_{1}}}\,\bar{b}\,\bar{a_{_{2}}}\,
\bar{b_{_{1}}}(c_{_{1,2}}) & \\

&=& c_{_{1,2}}\,b_{_{1}}\,a_{_{2}}\,b_{_{1}}\,c_{_{2,4}}\,
b_{_{1}}\,b_{_{2}}\,c_{_{2,4}}\,\bar{a_{_{1}}}\,\bar{b}\,\bar{a_{_{2}}}\,
\bar{b_{_{1}}}(c_{_{1,2}}) & \hbox{by {\it (T)}} \\

&=& c_{_{1,2}}\,b_{_{1}}\,a_{_{2}}\,c_{_{2,4}}\,b_{_{1}}\,c_{_{2,4}}\,
b_{_{2}}\,c_{_{2,4}}\,\bar{a_{_{1}}}\,\bar{b}\,\bar{a_{_{2}}}\,
\bar{b_{_{1}}}(c_{_{1,2}}) & \hbox{by {\it (T)}} \\

&=& c_{_{1,2}}\,b_{_{1}}\,a_{_{2}}\,c_{_{2,4}}\,b_{_{1}}\,b_{_{2}}\,
c_{_{2,4}}\,b_{_{2}}\,\bar{a_{_{1}}}\,\bar{b}\,\bar{a_{_{2}}}\,
\bar{b_{_{1}}}(c_{_{1,2}}) & \hbox{by {\it (T)}} \\

&=& c_{_{1,2}}\,b_{_{1}}\,a_{_{2}}\,c_{_{2,4}}\,b_{_{1}}\,b_{_{2}}\,
c_{_{2,4}}\,\bar{a_{_{1}}}\,\bar{b}\,\bar{a_{_{2}}}\,
\bar{b_{_{1}}}(c_{_{1,2}}) & \hbox{by {\it (T)}} \\

&=& x. & 
\end{array}$$

\noindent
To conclude, we get,
$$\begin{array}{rcll}
f & = & b_{_{2}}\,c_{_{2,4}}\,b_{_{1}}\,a_{_{2}}\,b\,\bar{c_{_{2,4}}}\,
\bar{b_{_{2}}}\,\bar{b_{_{1}}}\,\bar{a_{_{2}}}\,\bar{c_{_{2,4}}}\,
\bar{b_{_{1}}}(x) & \\

&=& b_{_{2}}\,c_{_{2,4}}\,b_{_{1}}\,a_{_{2}}\,b(x) & \\

&=& b_{_{2}}\,c_{_{2,4}}\,b_{_{1}}\,a_{_{2}}\,b\,c_{_{1,2}}\,b_{_{1}}\,
c_{_{2,4}}\,b_{_{2}}\,\bar{b}(a_{_{4}}) & \\

&=& b_{_{2}}\,c_{_{2,4}}\,b_{_{1}}\,a_{_{2}}\,c_{_{1,2}}\,b_{_{1}}\,
c_{_{2,4}}\,\bar{a_{_{4}}}(b_{_{2}}) & \hbox{by {\it (T)}} \\

&=& b_{_{2}}\,c_{_{2,4}}\,\bar{a_{_{4}}}\,\bar{b_{_{2}}}\,b_{_{1}}\,
a_{_{2}}\,c_{_{1,2}}\,b_{_{1}}(c_{_{2,4}}) & \hbox{by {\it (T)}} \\

&=& b_{_{2}}\,(a_{_{1}}\,a_{_{2}}\,a_{_{4}}\,b)^{3}\,\bar{c_{_{1,2}}}\,
\bar{c_{_{4,1}}}\,\bar{a_{_{4}}}\,\bar{b_{_{2}}}\,b_{_{1}}\,a_{_{2}}\,
c_{_{1,2}}\,b_{_{1}}(c_{_{2,4}}) & \hbox{by }\,(E_{_{1,2,4}}) \\

&=& b_{_{2}}\,(a_{_{1}}\,a_{_{2}}\,a_{_{4}}\,b)^{3}\,\bar{a_{_{4}}}\,
\bar{c_{_{4,1}}}\,\bar{b_{_{2}}}\,\bar{c_{_{1,2}}}\,b_{_{1}}\,
c_{_{1,2}}\,a_{_{2}}\,b_{_{1}}(c_{_{2,4}}) & \hbox{by {\it (T)}} \\

&=& b_{_{2}}\,(a_{_{1}}\,a_{_{2}}\,b)^{2}\,a_{_{4}}\,b\,a_{_{1}}\,a_{_{2}}\,b\,
\bar{c_{_{4,1}}}\,\bar{b_{_{2}}}\,\bar{c_{_{1,2}}}\,b_{_{1}}\,
c_{_{1,2}}\,a_{_{2}}\,b_{_{1}}(c_{_{2,4}}) & \hbox{by lemma~\ref{etoile}}\\

&=& (a_{_{1}}\,a_{_{2}}\,b)^{2}\,b_{_{2}}\,a_{_{4}}\,\bar{c_{_{4,1}}}\,
\bar{b_{_{2}}}\,b\,a_{_{1}}\,a_{_{2}}\,b\,
b_{_{1}}\,c_{_{1,2}}\,\bar{b_{_{1}}}\,a_{_{2}}\,
b_{_{1}}(c_{_{2,4}}) & \hbox{by {\it (T)}} \\

&=& (a_{_{1}}\,a_{_{2}}\,b)^{2}\,b_{_{2}}\,a_{_{4}}\,\bar{c_{_{4,1}}}\,
\bar{b_{_{2}}}\,b\,a_{_{2}}\,a_{_{1}}\,b\,
b_{_{1}}\,c_{_{1,2}}\,a_{_{2}}\,
b_{_{1}}\,\bar{a_{_{2}}}(c_{_{2,4}}) & \hbox{by {\it (T)}} \\

&=& (a_{_{1}}\,a_{_{2}}\,b)^{2}\,h(c_{_{2,4}}) &  \\

&=& (a_{_{1}}\,a_{_{2}}\,b)^{2}(c_{_{3,4}}) &  \\

&=& c_{_{3,4}} & \hbox{by {\it (T)}\,.}  \\
\end{array}$$
Finally, we have proved that $\,h(c_{_{2,4}})=c_{_{3,4}}=f$.

\vskip3mm\noindent
4) We can compute $\,h(a_{_{4}})\,$ as follows:
$$\begin{array}{rcll}
h(a_{_{4}}) & = & b_{_{2}}\,a_{_{4}}\,\bar{c_{_{4,1}}}\,\bar{b_{_{2}}}\,b
\,a_{_{2}}\,a_{_{1}}\,b\,b_{_{1}}\,c_{_{1,2}}\,a_{_{2}}\,b_{_{1}}(a_{_{4}}) &\\

&=& b_{_{2}}\,a_{_{4}}\,\bar{c_{_{4,1}}}\,\bar{b_{_{2}}}\,b
\,a_{_{2}}\,a_{_{1}}\,b(a_{_{4}}) & \hbox{by {\it (T)}} \\

&=& b_{_{2}}\,a_{_{4}}\,(\bar{a_{_{1}}}\,\bar{a_{_{2}}}\,\bar{a_{_{4}}}\,
\bar{b})^{3}\,c_{_{1,2}}\,c_{_{2,4}}\,\bar{b_{_{2}}}\,b
\,a_{_{2}}\,a_{_{1}}\,b(a_{_{4}}) & \hbox{by }\,(E_{_{1,2,4}}) \\

&=& b_{_{2}}\,c_{_{2,4}}\,\bar{a_{_{1}}}\,\bar{a_{_{2}}}\,
\bar{b}\,\bar{a_{_{1}}}\,\bar{a_{_{2}}}\,\bar{a_{_{4}}}\,
\bar{b}\,\bar{a_{_{1}}}\,\bar{a_{_{2}}}\,\bar{a_{_{4}}}\,
\bar{b}\,\bar{b_{_{2}}}\,b\,a_{_{2}}\,
a_{_{1}}\,b(a_{_{4}}) & \hbox{by {\it (T)}} \\

&=& b_{_{2}}\,c_{_{2,4}}\,\bar{a_{_{1}}}\,\bar{a_{_{2}}}\,
\bar{b}\,\bar{a_{_{1}}}\,\bar{a_{_{2}}}\,\bar{b}\,\bar{a_{_{4}}}\,
\bar{b}\,\bar{b_{_{2}}}\,b(a_{_{4}}) & \hbox{by {\it (T)}} \\

&=& b_{_{2}}\,c_{_{2,4}}\,\bar{a_{_{1}}}\,\bar{a_{_{2}}}\,
\bar{b}\,\bar{a_{_{1}}}\,\bar{a_{_{2}}}\,\bar{b}\,\bar{a_{_{4}}}\,
a_{_{4}}(b_{_{2}}) & \hbox{by {\it (T)}} \\

&=& b_{_{2}}\,c_{_{2,4}}(b_{_{2}}) & \hbox{by {\it (T)}} \\

&=& c_{_{2,4}} & \hbox{by {\it (T)}.} \\
\end{array}$$

\noindent
5) For $\,h(c_{_{1,4}})$, we have:
$$\begin{array}{rcll}
h(c_{_{1,4}}) & = & b_{_{2}}\,a_{_{4}}\,\bar{c_{_{4,1}}}\,\bar{b_{_{2}}}\,b\,
a_{_{2}}\,a_{_{1}}\,b\,b_{_{1}}\,c_{_{1,2}}\,a_{_{2}}\,b_{_{1}}(c_{_{1,4}}) & \\

&=& b_{_{2}}\,a_{_{4}}\,\bar{c_{_{4,1}}}\,b\,a_{_{2}}\,a_{_{1}}\,
b\,b_{_{1}}\,a_{_{2}}\,\bar{b_{_{2}}}(c_{_{1,4}}) & \hbox{by {\it (T)}} \\

&=& b_{_{2}}\,a_{_{4}}\,\bar{a_{_{4}}}\,\bar{b}\,\bar{a_{_{2}}}\,
\bar{a_{_{1}}}\,\bar{b}\,\bar{a_{_{4}}}\,\bar{a_{_{2}}}\,\bar{a_{_{1}}}\,
c_{_{1,2}}\,c_{_{2,4}}\,b_{_{1}}\,a_{_{2}}\,\bar{b_{_{2}}}(c_{_{1,4}})
& \hbox{by }\,(E_{_{1,2,4}}) \\

&=& \bar{b}\,\bar{a_{_{2}}}\,\bar{a_{_{1}}}\,\bar{b}\,\bar{a_{_{2}}}\,
c_{_{1,2}}\,b_{_{2}}\,c_{_{2,4}}\,b_{_{1}}\,\bar{a_{_{4}}}\,\bar{a_{_{1}}}\,
a_{_{2}}\,c_{_{1,4}}(b_{_{2}}) & \hbox{by {\it (T)}} \\

&=& \bar{b}\,\bar{a_{_{2}}}\,\bar{a_{_{1}}}\,\bar{b}\,\bar{a_{_{2}}}\,
c_{_{1,2}}\,b_{_{2}}\,c_{_{2,4}}\,b_{_{1}}\,c_{_{1,2}}\,c_{_{2,4}}\,\bar{X}\,
\bar{a_{_{2}}}\,X(b_{_{2}}) & \hbox{by }\,(L_{_{1,2,4}}) \\

&=& \bar{b}\,\bar{a_{_{2}}}\,\bar{a_{_{1}}}\,\bar{b}\,\bar{a_{_{2}}}\,
c_{_{1,2}}\,b_{_{2}}\,c_{_{2,4}}\,b_{_{1}}\,c_{_{2,4}}\,\bar{b}\,
\bar{a_{_{1}}}\,\bar{a_{_{4}}}\,\bar{b}\,\bar{a_{_{2}}}\,b\,a_{_{4}}(b_{_{2}})
& \hbox{by {\it (T)}} \\

&=& \bar{b}\,\bar{a_{_{2}}}\,\bar{a_{_{1}}}\,\bar{b}\,\bar{a_{_{2}}}\,
c_{_{1,2}}\,b_{_{2}}\,b_{_{1}}\,c_{_{2,4}}\,b_{_{1}}\,\bar{b}\,
\bar{a_{_{1}}}\,\bar{a_{_{4}}}\,a_{_{2}}\,\bar{b}\,\bar{a_{_{2}}}\,
a_{_{4}}(b_{_{2}}) & \hbox{by {\it (T)}} \\

&=& \bar{b}\,\bar{a_{_{2}}}\,\bar{a_{_{1}}}\,\bar{b}\,\bar{a_{_{2}}}\,
c_{_{1,2}}\,b_{_{1}}\,b_{_{2}}\,c_{_{2,4}}\,b_{_{1}}\,\bar{b}\,
\bar{a_{_{1}}}\,a_{_{2}}\,b\,\bar{a_{_{4}}}\,\bar{b}(b_{_{2}})
& \hbox{by {\it (T)}} \\

&=& \bar{b}\,\bar{a_{_{2}}}\,\bar{a_{_{1}}}\,\bar{b}\,\bar{a_{_{2}}}\,
c_{_{1,2}}\,b_{_{1}}\,b_{_{2}}\,c_{_{2,4}}\,b_{_{1}}\,\bar{b}\,
\bar{a_{_{1}}}\,a_{_{2}}\,b\,b_{_{2}}(a_{_{4}}) & \hbox{by {\it (T)}} \\

&=& \bar{b}\,\bar{a_{_{2}}}\,\bar{a_{_{1}}}\,\bar{b}\,\bar{a_{_{2}}}\,
c_{_{1,2}}\,b_{_{1}}\,c_{_{2,4}}\,b_{_{2}}\,c_{_{2,4}}\,b_{_{1}}\,\bar{b}\,
\bar{a_{_{1}}}\,a_{_{2}}\,\bar{a_{_{4}}}(b) & \hbox{by {\it (T)}} \\

&=& \bar{b}\,\bar{a_{_{2}}}\,\bar{a_{_{1}}}\,\bar{b}\,\bar{a_{_{2}}}\,
c_{_{1,2}}\,b_{_{1}}\,c_{_{2,4}}\,b_{_{2}}\,c_{_{2,4}}\,\bar{b}\,\bar{a_{_{1}}}\,
\bar{a_{_{4}}}\,\bar{b}\,b_{_{1}}(a_{_{2}}) & \hbox{by {\it (T)}} \\

&=& \bar{b}\,\bar{a_{_{2}}}\,\bar{a_{_{1}}}\,\bar{b}\,\bar{a_{_{2}}}\,
c_{_{1,2}}\,b_{_{1}}\,c_{_{2,4}}\,b_{_{2}}\,c_{_{2,4}}\,\bar{b}\,\bar{a_{_{1}}}\,
\bar{a_{_{4}}}\,\bar{b}\,\bar{a_{_{2}}}(b_{_{1}}) & \hbox{by {\it (T)}} \\

&=& \bar{b}\,\bar{a_{_{2}}}\,\bar{a_{_{1}}}\,\bar{b}\,\bar{a_{_{2}}}\,
c_{_{1,2}}\,b_{_{1}}\,c_{_{2,4}}\,b_{_{2}}\,a_{_{1}}\,a_{_{4}}\,a_{_{2}}\,X\,
\bar{c_{_{1,2}}}\,\bar{c_{_{4,1}}}(b_{_{1}})
& \hbox{by }\,(E_{_{1,2,4}}) \\

&=& \bar{b}\,\bar{a_{_{2}}}\,\bar{a_{_{1}}}\,\bar{b}\,c_{_{1,2}}\,
\bar{a_{_{2}}}\,b_{_{1}}\,a_{_{2}}\,\bar{c_{_{1,2}}}\,c_{_{2,4}}(b_{_{1}})
& \hbox{by {\it (T)}} \\

&=& \bar{b}\,\bar{a_{_{2}}}\,\bar{a_{_{1}}}\,\bar{b}\,c_{_{1,2}}\,
b_{_{1}}\,a_{_{2}}\,\bar{b_{_{1}}}\,\bar{c_{_{1,2}}}\,\bar{b_{_{1}}}(c_{_{2,4}})
& \hbox{by {\it (T)}} \\

&=& \bar{b}\,\bar{a_{_{2}}}\,\bar{a_{_{1}}}\,\bar{b}\,c_{_{1,2}}\,
b_{_{1}}\,a_{_{2}}\,\bar{c_{_{1,2}}}\,\bar{b_{_{1}}}\,
\bar{c_{_{1,2}}}(c_{_{2,4}}) & \hbox{by {\it (T)}} \\

&=& \bar{b}\,\bar{a_{_{2}}}\,\bar{a_{_{1}}}\,\bar{b}\,\bar{b_{_{1}}}\,
c_{_{1,2}}\,b_{_{1}}\,a_{_{2}}\,\bar{b_{_{1}}}(c_{_{2,4}})
& \hbox{by {\it (T)}} \\

&=& \bar{b}\,\bar{a_{_{2}}}\,\bar{a_{_{1}}}\,\bar{b}\,\bar{b_{_{1}}}\,
c_{_{1,2}}\,\bar{a_{_{2}}}\,b_{_{1}}\,a_{_{2}}(c_{_{2,4}})
& \hbox{by {\it (T)}} \\

&=& \bar{b}\,\bar{a_{_{2}}}\,\bar{a_{_{1}}}\,\bar{b}\,\bar{b_{_{1}}}\,
\bar{a_{_{2}}}\,\bar{c_{_{2,4}}}\,c_{_{1,2}}(b_{_{1}})
& \hbox{by {\it (T)}} \\

&=& \bar{b}\,\bar{a_{_{2}}}\,\bar{a_{_{1}}}\,\bar{b}\,\bar{b_{_{1}}}\,
\bar{a_{_{2}}}\,\bar{c_{_{2,4}}}\,\bar{b_{_{1}}}(c_{_{1,2}})
& \hbox{by {\it (T)}} \\

&=& l\,.
\end{array}$$

\noindent
6) By the relations {\it (T)}, one has
$$\begin{array}{rcl}
h(a_{_{2}}) & = & b_{_{2}}\,a_{_{4}}\,\bar{c_{_{4,1}}}\,\bar{b_{_{2}}}\,b
\,a_{_{2}}\,a_{_{1}}\,b\,b_{_{1}}\,c_{_{1,2}}\,a_{_{2}}\,b_{_{1}}(a_{_{2}}) \\

&=& b_{_{2}}\,a_{_{4}}\,\bar{c_{_{4,1}}}\,\bar{b_{_{2}}}\,b\,a_{_{2}}\,a_{_{1}}\,
b\,b_{_{1}}\,c_{_{1,2}}\,a_{_{2}}\,\bar{a_{_{2}}}(b_{_{1}}) \\

&=& b_{_{2}}\,a_{_{4}}\,\bar{c_{_{4,1}}}\,\bar{b_{_{2}}}\,b
\,a_{_{2}}\,a_{_{1}}\,b\,b_{_{1}}\,\bar{b_{_{1}}}(c_{_{1,2}}) \\

&=& c_{_{1,2}}\,.
\end{array}$$

\noindent
7) Using the braid relations, one gets
$$\begin{array}{rcl}
h(b) & = & b_{_{2}}\,a_{_{4}}\,\bar{c_{_{4,1}}}\,\bar{b_{_{2}}}\,b
\,a_{_{2}}\,a_{_{1}}\,b\,b_{_{1}}\,c_{_{1,2}}\,a_{_{2}}\,b_{_{1}}(b) \\

&=& b_{_{2}}\,a_{_{4}}\,\bar{c_{_{4,1}}}\,\bar{b_{_{2}}}\,b
\,a_{_{2}}\,a_{_{1}}\,b\,b_{_{1}}\,\bar{b}(a_{_{2}}) \\

&=& b_{_{2}}\,a_{_{4}}\,\bar{c_{_{4,1}}}\,\bar{b_{_{2}}}\,b
\,a_{_{2}}\,\bar{a_{_{2}}}(b_{_{1}}) \\

&=& b_{_{1}}\,.\\
\end{array}$$

\vskip3mm\noindent
Thus, one has 
$\,h\bar{X}(a_{_{2}})=\bar{b_{_{1}}}\,\bar{a_{_{2}}}\,\bar{c_{_{2,4}}}\,
\bar{b_{_{1}}}(c_{_{1,2}})=m$.

\vskip3mm\noindent
This concludes the proof of lemma~\ref{psi}.

\eproof


\section{Proof of theorem~\ref{principaltheorem} \label{finpreuve}}

We will proceed by induction on $n$. To do this, we need the exact 
sequence (see \cite{Birman,Harer}):
\diagram[size=1.5em]
 1 & \rto &   \mathbf{Z}\times\pi_{1}(\Sigma_{g,n-1},p)  & \rto^{f_{1}} &
  \mathcal{M}_{g,n} & \rto^{f_{2}} & \mathcal{M}_{g,n-1} & \rto & 1\ .  \\
\enddiagram
Here, $f_{2}\,$ is defined by collapsing $\delta_{n}\,$ with a disc 
centred at $p$ and by extending each map over the disc by the identity, and 
$f_{1}\,$ by sending each $\,k\in\mathbf{Z}\,$ to 
$\,\tau_{\delta_{n}}^k\,$ and each $\,\alpha\in\pi_{1}(\Sigma_{g,n-1},p)\,$ 
to the spin map $\,\tau_{\alpha'}\tau_{\alpha''}^{-1}\,$ $\,$($\alpha'$ and $\alpha''$ 
are two curves in $\Sigma_{g,n-1}\,$ which are separated by $\delta_{n}\,$ 
and such that $\,\alpha'=\alpha''=\alpha\,$ in $\Sigma_{g,n-1}$).

Let us denote by $\,a'_{_{1}},\ldots,a'_{_{2g+n-3}},b',b'_{_{1}},\ldots,
b'_{_{g-1}},(c'_{_{i,j}})_{1\leq i\not= j\leq2g+n-3}\,$ the \linebreak[4]
generators of $\,G_{g,n-1}\,$ corresponding to the curves in
$\,\mathcal{G}_{g,n-1}$. We define $\,g_{2}:G_{g,n}\rightarrow G_{g,n-1}\,$ by

$$
\begin{array}{rcll}
g_{2}(a_{_{i}}) & = & a'_{_{i}} & \hbox{ for all }\,  i\not= 2g+n-2 \\
g_{2}(a_{_{2g+n-2}}) & = & a'_{_{1}} &    \\
g_{2}(b) & = & b'   & \\
g_{2}(b_{_{i}}) & = & b'_{_{i}} & \hbox{ for }\, 1\leq i\leq g-1 \\
g_{2}(c_{_{i,j}}) & = & c'_{_{i,j}}& \hbox{ for }\, 1\leq i,j\leq 2g+n-3 \\
g_{2}(c_{_{i,2g+n-2}}) & = & c'_{_{i,1}}& \hbox{ for }\, 2\leq i\leq 2g+n-3 \\
g_{2}(c_{_{2g+n-2,j}}) & = & c'_{_{1,j}}&\hbox{ for }\, 2\leq j\leq 2g+n-3 \\
g_{2}(c_{_{1,2g+n-2}}) & = & (a'_{_{1}}\,b'\,a'_{_{1}})^{4} \\
g_{2}(c_{_{2g+n-2,1}}) & = & 1\,.&
\end{array}
$$

\vskip3mm
\begin{lemma}
For all $\,(g,n)\!\in\!{\mathbf{N}}^{\ast}\!\times\!{\mathbf{N}}^{\ast}$,
$g_{2}\,$ is an homomorphism.
\end{lemma}

\proof We have to prove that the relations in $\,G_{g,n}\,$ are 
satisfied in $\,G_{g,n-1}\,$ via $g_{2}$. Since for all $i$ such 
that $\,1\leq i\leq g-1$, one has $\,g_{2}(c_{_{2i,2i+1}})\!=\!c'_{_{2i,2i+1}}\,$ 
and $\,g_{2}(c_{_{2i-1,2i}})\!=\!c'_{_{2i-1,2i}}\,$, this is clear for the 
handle relations.

So, let $\lambda$, $\mu$ be two elements of $\,\mathcal{G}_{g,n}\,$ which 
do not intersect (resp. intersect transversaly in a single point). 
If $\,l\,$ and $\,m\,$ are the associated elements of $\,G_{g,n}$, we have 
to prove that

\[
(\bullet)\left\{\begin{array}{c}
\,g_{2}(l)g_{2}(m)=g_{2}(m)g_{2}(l) \\
\Bigl(\hbox{resp. }\ g_{2}(l)g_{2}(m)g_{2}(l)=g_{2}(m)g_{2}(l)g_{2}(m)\Bigr).
\end{array}\right.\]

\vskip3mm\noindent
When $\lambda$ and $\mu$ are distinct from $\,\gamma_{_{2g+n-2,1}}\,$ and 
$\,\gamma_{_{1,2g+n-2}}$, these relations are precisely braid relations 
in $\,G_{g,n-1}$. If not, $\lambda$ and $\mu$ do not intersect in a 
single point. Thus, it remains to consider the cases where $\,\lambda=\gamma_{_{1,2g+n-2}}\,$
or $\,\gamma_{_{2g+n-2,1}}\,$ and $\,\mu\in\mathcal{G}_{g,n}\,$ is a curve 
disjoint from $\lambda$. For $\,\lambda=\gamma_{_{2g+n-2,1}}$, one has 
$\,g_{2}(l)=1\,$ and the relation $\,(\bullet)\,$ is satisfied in 
$\,G_{g,n-1}$. So, suppose that $\,\lambda\!=\!\gamma_{_{1,2g+n-2}}$. Then, 
we have $\,g_{2}(l)\!=\!(a'_{_{1}}\,b'\,a'_{_{1}})^{4}$. The 
curves in $\,\mathcal{G}_{g,n}\,$ which are disjoint from $\lambda$ 
are 
$\,\beta,\beta_{1},\ldots,\beta_{g-1},\alpha_{1},\alpha_{2g+n-2},\gamma_{2g+n-2,1}\,$
and $\,(\gamma_{i,j})_{1\leq i<j\leq 2g+n-2}$. Let us look at the 
different cases:
\begin{list}{--}{\leftmargin10mm\parsep=2mm\topsep=4mm}
\item By lemma~\ref{etoile}, $\,b'=g_{2}(b)\,$ and $\,a'_{_{1}}=
      g_{2}(a_{_{1}})=g_{2}(a_{_{2g+n-2}})\,$ commute with $\,
      (a'_{_{1}}\,b'\,a'_{_{1}})^{4}=g_{2}(l)$.
\item For all $i$, $\,1\leq i\leq g-1$, 
      $\,b'_{_{i}}=g_{2}(b_{_{i}})\,$ commutes with $\,
      (a'_{_{1}}\,b'\,a'_{_{1}})^{4}$ by the braid relations 
      in $\,G_{g,n-1}$.
\item For all $\,i,j\,$ such that $\,1\leq i<j\leq 2g+n-2$, one has 
$\,g_{2}(c_{_{i,j}})=c'_{_{i,j}}\,$ if $\,j\not = 2g+n-2$, and
$\,g_{2}(c_{_{i,j}})=c'_{_{i,1}}\,$ otherwise. In all cases, one 
has that $\,g_{2}(c_{_{i,j}})g_{2}(l)=g_{2}(l)g_{2}(c_{_{i,j}})\,$ by the 
braid relations in $\,G_{g,n-1}$.
\end{list}

\vskip3mm
Now, let us look at the star relations. For $\,i,j,k\!\not =\! 2g+n-2$,
$\,(E_{_{i,j,k}})$ is sent by $g_{2}\,$ to
$\,(E'_{_{i,j,k}})$, the star relation in $\,G_{g,n-1}\,$ involving
the same curves. For all $i,j$ such that $\,2\leq i\leq j<2g+n-2$, 
$\,(E_{_{i,j,2g+n-2}})\,$ is sent to $\,(E'_{_{i,j,1}})$. Next, for 
$\,2\leq j<2g+n-2$, $\,(E_{_{1,j,2g+n-2}})\,$ is sent to $\,(E'_{_{1,1,j}})$.
Finally, since $\,g_{2}(c_{_{2g+n-2,1}})=1\,$ and 
$\,g_{2}(c_{_{1,2g+n-2}})=(a'_{_{1}}b'a'_{_{1}})^{4}$, the 
relation $\,(E_{_{1,1,2g+n-2}})\,$ is satisfied in $\,G_{g,n-1}\,$ 
via $g_{2}\,$ by lemma~\ref{etoile}. This concludes the proof by remark~\ref{rem}.

\eproof

\hfill\break\indent
Since the relations {\it (T)}, {\it (A)} and {\it 
(E$_{i,j,k}$)} are satisfied in $\,\mathcal{M}_{g,n}\,$ (see \cite{Gervais}),
one has an homomorphism $\,\Phi_{g,n}:G_{g,n}\rightarrow \mathcal{M}_{g,n}\,$
which associates to each $\,a\in\mathcal{G}_{g,n}\,$ the corresponding twist
$\,\tau_{\alpha}$. Since we view $\,\Sigma_{g,n}\,$ as a subsurface of
$\,\Sigma_{g,n-1}$, we have $\,\Phi_{g,n-1}\circ g_{2}\!=\!f_{2}\circ\Phi_{g,n}$.
Thus, we get the following commutative diagram:

\diagram[size=2.5em]
 1 & \rto &   \ker g_{2}   & \rto & G_{g,n} &  \rTo^{g_{2}} &&
                                                        G_{g,n-1} & \rto & 1  \\                                                        
   &      & \dto_{h_{g,n}}  &   & \dto_{\Phi_{g,n}} & &&
                                                 \dto_{\Phi_{g,n-1}} & &\\
 1 & \rto & \mathbf{Z}\times\pi_{1}(\Sigma_{g,n-1},p) & \rTo^{f_{1}} &
         \mathcal{M}_{g,n} & \rto^{f_{2}} && \mathcal{M}_{g,n-1} & \rto & 1 \\
\enddiagram

\vskip3mm\noindent
where $h_{g,n}\,$ is induced by $\Phi_{g,n}$.

\begin{proposition}\label{h}
$h_{g,n}\,$ is an isomorphism for all $\,g\!\geq\!1$ and $\,n\!\geq\! 2$.
\end{proposition}

In order to prove this proposition, we will first give a system of 
generators for $\ker g_{2}$. Thus, we consider the following elements of
$\,\ker g_{2}$:
$$\begin{array}{c}
x_{_{0}}=a_{_{1}}\bar{a_{_{2g+n-2}}},\ \ x_{_{1}}=b(x_{_{0}}),\ \ 
x_{_{2}}=a_{_{2}}(x_{_{1}}),\ \ x_{_{3}}=b_{_{1}}(x_{_{2}}),\\ \\
\hbox{for }\, 2\leq i\leq g-1,\ 
\,x_{_{2i}}=c_{_{2i-2,2i}}(x_{_{2i-1}})\ \hbox{ and }\ 
x_{_{2i+1}}=b_{_{i}}(x_{_{2i}}),\\ \\
\hbox{and for }\,2g\leq k\leq 2g+n-3,\ \ 
x_{_{k}}=a_{_{k}}(x_{_{1}})\,.
\end{array}$$

\begin{remark}
If $\,g\!=\!1$, one has just to concider 
$\,x_{_{0}},x_{_{1}},x_{_{2}},\ldots,x_{_{n-1}}$.
\end{remark}

\begin{lemma}\label{normal}
For all $\,(g,n)\!\in\!{\mathbf{N}}^{\ast}\!\times\!{\mathbf{N}}^{\ast}$,
$\,\ker g_{2}\,$ is normally generated by $\,d_{_{n}}\,$ and $\,x_{_{0}}$.
\end{lemma}

\proof Let us denote by $K$ the subgroup of $\,G_{g,n}\,$ normally 
generated by $\,d_{_{n}}\,$ and $\,x_{_{0}}$. Since 
$\,g_{2}(d_{_{n}})\!=\! 1\,$ and $\,g_{2}(a_{_{2g+n-2}})\!=\!
g_{2}(a_{_{1}})$, one has $\,K\!\subset\!\ker g_{2}$. In order to 
prove the equality, we shall prove that $\,g_{2}\,$ induces a 
monomorphism $\,\widetilde{g_{2}}\,$ from $\,G_{g,n}/K\,$ to 
$\,G_{g,n-1}$.

\noindent
Define $\,k:G_{g,n-1}\rightarrow G_{g,n}/K\,$ by
$$
\begin{array}{rcll}
k(b') & = & \widetilde{b}   & \\
k(b'_{_{i}}) & = & \widetilde{b_{_{i}}} & \hbox{ for }\, 1\leq i\leq g-1 \\
k(a'_{_{i}}) & = & \widetilde{a_{_{i}}} & 
              \hbox{ for all }i,\,\  1\leq i\leq 2g+n-3 \\
k(c'_{_{i,j}}) & = & \widetilde{c_{_{i,j}}} & 
              \hbox{ for all }\,i\not= j,\,\ 1\leq i,j\leq 2g+n-3
\end{array}
$$
where, for $\,x\!\in\! G_{g,n}$, $\,\widetilde{x}$ denote the class of 
$x$ in $\,G_{g,n}/K$. Pasting a pair of pants to $\,\gamma_{2g+n-3,1}\,$ allows 
us to view $\,\Sigma_{g,n-1}\,$ as a subsurface of $\,\Sigma_{g,n}$, 
and $\,\mathcal{G}_{g,n-1}\,$ as a subset of $\,\mathcal{G}_{g,n}$. 
Thus, $k$ appears to be clearly a morphism. Let us prove that
$\,k\!\circ \widetilde{g_{2}}\!=\!Id$.

\noindent
Denote by $H$ the subgroup of $\,G_{g,n}/K\,$ generated by $\,\{\widetilde{b},
\widetilde{b_{_{1}}},\ldots,\widetilde{b}_{_{g-1}},$
$\widetilde{a_{_{1}}},\ldots,\widetilde{a}_{_{2g+n-3}},
(\widetilde{c}_{_{i,j}})_{_{1\leq i\not= j\leq 2g+n-3}}\}$. 
Since, by definition of $\,g_{2}\,$ and $k$, one has $\,k\!\circ\! 
g_{2}(\widetilde{x})\!=\!\widetilde{x}\,$ for all 
$\,\widetilde{x}\!\in\!H$, we just need to prove that\linebreak[4] 
$\,G_{g,n}/K=H$. We know that $\,G_{g,n}/K\,$ is generated by 
$\,\{\widetilde{x}\,/\,x\!\in\!\mathcal{G}_{g,n}\}$; thus, the 
following computations allow us to conclude.

\vskip3mm
\begin{list}{--}{\leftmargin10mm\itemsep3mm}
\item $\widetilde{a}_{_{2g+n-2}}\!=\,\widetilde{a_{_{1}}}$.
\item $\widetilde{c}_{_{2g+n-2,1}}\!=\!\widetilde{d_{_{n}}}\!=\!1$.
\item By the star relation $\,(E_{_{1,1,2g+n-2}})$, one has
      $$\widetilde{c}_{_{1,2g+n-2}}\!=\!(\widetilde{a_{_{1}}}\,
      \widetilde{a_{_{1}}}\,\widetilde{a}_{_{2g+n-2}}\,\widetilde{b})^{-3}\,
      \widetilde{c}_{_{2g+n-2,1}}\,=\,(\widetilde{a_{_{1}}}\,
      \widetilde{a_{_{1}}}\,\widetilde{a_{_{1}}}\,\widetilde{b})^{-3}\,.$$
\item For $\,2\!\leq\! i\!\leq\! 2g+n-3$, one has by the lantern relation
      $\,(L_{_{2g+n-2,1,i}})$:
      $$a_{_{2g+n-2}}\,c_{_{2g+n-2,1}}\,c_{_{1,i}}\,a_{_{i}}=c_{_{2g+n-2,i}}\,
      a_{_{1}}\,X\,a_{_{1}}\,\bar{X}$$
      where $\,X\!=\!b\,a_{_{2g+n-2}}\,a_{_{i}}\,b$.
      This relation implies the following one by $\,${\it (T)}:
      $$\begin{array}{rcl}
      c_{_{2g+n-2,i}} & = & c_{_{1,i}}\,a_{_{i}}\,X\,\bar{a_{_{1}}}\,\bar{X}\,
          \bar{a_{_{1}}}\,a_{_{2g+n-2}}\,c_{_{2g+n-2,1}}\\
      &=& c_{_{1,i}}\,X\,\bar{x_{_{0}}}\,\bar{X}\,\bar{x_{_{0}}}\,d_{_{n}}\,,
      \end{array}$$
      which yields $\,\widetilde{c}_{_{2g+n-2,i}}\!=\!\widetilde{c}_{_{1,i}}$.
\item In the same way, using the lantern relation $\,(L_{_{i,2g+n-2,1}})$,
      one proves that $\,\widetilde{c}_{_{i,2g+n-2}}\!=\!
      \widetilde{c}_{_{i,1}}\,$ for $\,2\!\leq\! i\!\leq\! 2g+n-3$.
\end{list}
\eproof

\begin{lemma}\label{gen-ker g2}
For all $\,(g,n)\!\in\!{\mathbf{N}}^{\ast}\!\times\!{\mathbf{N}}^{\ast}$,
$\,\ker g_{2}\,$ is generated by\linebreak[4] $\,d_{_{n}}\!=\!c_{_{2g+n-2,1}}\,$ and
$\,x_{_{0}},\ldots,x_{_{2g+n-3}}$.
\end{lemma}

\proof By lemma~\ref{normal}, $\ker g_{2}\,$ is normally generated by
$\,d_{_{n}}\,$ and $\,x_{_{0}}$. Furthermore, by the braid relations,
$\,d_{_{n}}\,$ is central in $\,G_{g,n}$. Thus, denoting by $K$ the subgroup
generated by $\,d_{_{n}},x_{_{0}},\ldots,x_{_{2g+n-2}}$, we have to 
prove that $\,gx_{_{0}}g^{-1}\!\in\! K\,$ for all $\,g\!\in\!G_{g,n}$. 
To do this, it is enough to show that $K$ is a normal subgroup of $\,G_{g,n}$.

By proposition~\ref{generator}, $\,G_{g,n}\,$ is 
generated by $\,\mathcal{H}_{g,n}\!=\!\{a_{_{1}},b,a_{_{2}},b_{_{1}},\ldots,$
$b_{_{g-1}},c_{_{2,4}},\ldots,c_{_{2g-4,2g-2}},c_{_{1,2}},a_{_{2g}},\ldots,
a_{_{2g+n-2}},d_{_{1}},\ldots,d_{_{n-1}}\}$. Since, by the braid relations, 
$\,d_{_{1}},\ldots,d_{_{n-1}}\,$ are central in $\,G_{g,n}$, we 
have to prove that $\,y(x_{_{k}})\,$ and $\,\bar{y}(x_{_{k}})\,$ are 
elements of $K$ for all $k$, $\,0\leq k\leq 
2g+n-3$, and all $\,y\!\in\! \mathcal{E}\,$ where 
$\,\mathcal{E}=\mathcal{H}_{g,n}\!\setminus\!\{d_{_{1}},\ldots,d_{_{n-1}}\}$.

\vskip3mm\noindent
$\ast$ \underline{Case 1: $\,k\!=\!0$}.

\vskip3mm\begin{list}{--}{\leftmargin7mm\itemsep3mm}
\item $\,b(x_{_{0}})=x_{_{1}}$.
\item We prove, using relations {\it (T)}, that
      $\,\bar{b}(x_{_{0}})=x_{_{0}}\,\bar{x_{_{1}}}\,x_{_{0}}$:
      $$\begin{array}{rcl}
      x_{_{0}}\,\bar{x_{_{1}}}\,x_{_{0}} & = & 
      a_{_{1}}\,\bar{a_{_{2g+n-2}}}\,b\,a_{_{2g+n-2}}\,
      \bar{a_{_{1}}}\,\bar{b}\,a_{_{1}}\,\bar{a_{_{2g+n-2}}} \\
      
      &=& a_{_{1}}\,b\,a_{_{2g+n-2}}\,\bar{b}\,b\,
      \bar{a_{_{1}}}\,\bar{b}\,\bar{a_{_{2g+n-2}}} \\
      
      &=& \bar{b}\,a_{_{1}}\,b\,\bar{b}\,\bar{a_{_{2g+n-2}}}\,b  \\
      
      &=& \bar{b}(x_{_{0}})\,.
      \end{array}$$
\item For $\,y\!\in\!\mathcal{E}\!\setminus\!\{b\}$, one has 
      $\,y(x_{_{0}})\!=\bar{y}(x_{_{0}})\!=\!x_{_{0}}\,$ by the braid
      relations.
\end{list}

\vskip5mm\noindent
$\ast$ \underline{Case 2: $\,k\!=\!1$}.

\vskip3mm\begin{list}{--}{\leftmargin7mm\itemsep3mm}
\item $a_{_{1}}(x_{_{1}})$ \parbox[t]{73mm}{$=a_{_{1}}\,b\,a_{_{1}}\,
      \bar{a_{_{2g+n-2}}}\,\bar{b}\,\bar{a_{_{1}}}=b\,a_{_{1}}\,b\,
      \bar{a_{_{2g+n-2}}}\,\bar{b}\,\bar{a_{_{1}}}$ $=b\,a_{_{1}}\,
      \bar{a_{_{2g+n-2}}}\,\bar{b}\,a_{_{2g+n-2}}\,\bar{a_{_{1}}}=x_{_{1}}\,
      \bar{x_{_{0}}}\,$,}
\item[] $\bar{a_{_{1}}}(x_{_{1}})$ \parbox[t]{73mm}{$=\bar{a_{_{1}}}\,b\,a_{_{1}}\,
      \bar{a_{_{2g+n-2}}}\,\bar{b}\,a_{_{1}}=b\,a_{_{1}}\,\bar{b}\,
      \bar{a_{_{2g+n-2}}}\,\bar{b}\,a_{_{1}}$ $=b\,a_{_{1}}\,
      \bar{a_{_{2g+n-2}}}\,\bar{b}\,\bar{a_{_{2g+n-2}}}\,a_{_{1}}=x_{_{1}}\,
      x_{_{0}}\,$.}
\item $a_{_{2g+n-2}}(x_{_{1}})$ \parbox[t]{60mm}{$=a_{_{2g+n-2}}\,b\,a_{_{1}}\,
      \bar{a_{_{2g+n-2}}}\,\bar{b}\,\bar{a_{_{2g+n-2}}}$
      
      $=a_{_{2g+n-2}}\,b\,a_{_{1}}\,\bar{b}\,\bar{a_{_{2g+n-2}}}\,\bar{b}$ 
      
      $=a_{_{2g+n-2}}\,\bar{a_{_{1}}}\,
      b\,a_{_{1}}\,\bar{a_{_{2g+n-2}}}\,\bar{b}=\bar{x_{_{0}}}\,x_{_{1}}\,$,}
\item[] $\bar{a_{_{2g+n-2}}}(x_{_{1}})$ \parbox[t]{60mm}{$=\bar{a_{_{2g+n-2}}}\,
      b\,a_{_{1}}\,\bar{a_{_{2g+n-2}}}\,\bar{b}\,a_{_{2g+n-2}}$
      
      $=\bar{a_{_{2g+n-2}}}\,
      b\,a_{_{1}}\,b\,\bar{a_{_{2g+n-2}}}\,\bar{b}$
      
      $=\bar{a_{_{2g+n-2}}}\,a_{_{1}}\,
      b\,a_{_{1}}\,\bar{a_{_{2g+n-2}}}\,\bar{b}=x_{_{0}}\,x_{_{1}}\,$.}
\item One has $\bar{b}(x_{_{1}})\!=\!x_{_{0}}\,$, and by the braid relations,
      $\,b(x_{_{1}})\!=\!x_{_{1}}\,\bar{x_{_{0}}}\,x_{_{1}}\,$:
      $$\begin{array}{rcl}
          x_{_{1}}\,\bar{x_{_{0}}}\,x_{_{1}} & = & b\,a_{_{1}}\,
          \bar{a_{_{2g+n-2}}}\,\bar{b}\,\bar{a_{_{1}}}\,a_{_{2g+n-2}}\,b\,
          a_{_{1}}\,\bar{a_{_{2g+n-2}}}\,\bar{b} \\
          
          &=& b\,\bar{a_{_{2g+n-2}}}\,\bar{b}\,\bar{a_{_{1}}}\,b\,\bar{b}\,
          a_{_{2g+n-2}}\,b\,a_{_{1}}\,\bar{b} \\
          
          &=& b\,b\,\bar{a_{_{2g+n-2}}}\,\bar{b}\,b\,a_{_{1}}\,\bar{b}\,\bar{b} \\
          
          &=& b(x_{_{1}}).
      \end{array}$$
\item For $\,i\!\in\!\{2,2g,2g+1,\ldots,2g+n-3\}$, we have 
      $\,a_{_{i}}(x_{_{1}})\!=\!x_{_{i}}\,$ and $\,\bar{a_{_{i}}}(x_{_{1}})\!=
      \!x_{_{1}}\,\bar{x_{_{i}}}\,x_{_{1}}\,$:
      
      $$\begin{array}{rcll}
      x_{_{1}}\,\bar{x_{_{i}}}\,x_{_{1}} & = & 
      b\,x_{_{0}}\,\bar{b}\,a_{_{i}}\,b\,\bar{x_{_{0}}}\, 
      \bar{b}\,\bar{a_{_{i}}}\,b\,x_{_{0}}\,\bar{b} & \\ 

      &=& b\,x_{_{0}}\,a_{_{i}}\,b\,\bar{a_{_{i}}}\,\bar{x_{_{0}}}\, 
      a_{_{i}}\,\bar{b}\,\bar{a_{_{i}}}\,x_{_{0}}\,\bar{b} & 
      \hbox{by {\it (T)}} \\ 

      &=& b\,a_{_{i}}\,x_{_{0}}\,b\,\bar{x_{_{0}}}\,\bar{b}\,x_{_{0}}\,
      \bar{a_{_{i}}}\,\bar{b} & \hbox{by case 1} \\ 

      &=& b\,a_{_{i}}\,x_{_{0}}\,\bar{x_{_{1}}}\,x_{_{0}}\,
      \bar{a_{_{i}}}\,\bar{b} & \\ 

      &=& b\,a_{_{i}}\,\bar{b}\,x_{_{0}}\,b\,
      \bar{a_{_{i}}}\,\bar{b} & \hbox{by case 1} \\ 

      &=& \bar{a_{_{i}}}\,b\,a_{_{i}}\,x_{_{0}}\,
      \bar{a_{_{i}}}\,\bar{b}\,a_{_{i}} & \hbox{by {\it (T)}} \\ 

      &=& \bar{a_{_{i}}}(x_{_{1}}) & \hbox{by case 1}.
      \end{array}$$
\item Each $\,y\!\in\!\{b_{_{1}},\ldots,b_{_{g-1}},c_{_{2,4}},\ldots,
      c_{_{2g-4,2g-2}},c_{_{1,2}}\}\,$ commutes with $x_{_{1}}\,$ by the 
      braid relations, so $\,y(x_{_{1}})\!=\bar{y}(x_{_{1}})\!=\!x_{_{1}}\,$.
\end{list}

\vskip5mm\noindent
$\ast$ \underline{Case 3: $\,k\!\in\!\{2,2g,\ldots,2g+n-3\}$}.

\vskip3mm\begin{list}{--}{\leftmargin7mm\itemsep3mm}
\item By the braid relations and the preceeding cases, we have:
      $$a_{_{1}}(x_{_{k}})=a_{_{k}}\,a_{_{1}}(x_{_{1}})=a_{_{k}}\,
      x_{_{1}}\,\bar{x_{_{0}}}\,\bar{a_{_{k}}}=x_{_{k}}\,\bar{x_{_{0}}}\,,$$
      $$\bar{a_{_{1}}}(x_{_{k}})=a_{_{k}}\,\bar{a_{_{1}}}(x_{_{1}})=a_{_{k}}\,
      x_{_{1}}\,x_{_{0}}\,\bar{a_{_{k}}}=x_{_{k}}\,x_{_{0}}\,,$$
      $$a_{_{2g+n-2}}(x_{_{k}})=a_{_{k}}\,a_{_{2g+n-2}}(x_{_{1}})=a_{_{k}}\,
      \bar{x_{_{0}}}\,x_{_{1}}\,\bar{a_{_{k}}}=\bar{x_{_{0}}}\,x_{_{k}}\,,$$
      $$\bar{a_{_{2g+n-2}}}(x_{_{k}})=a_{_{k}}\,\bar{a_{_{2g+n-2}}}(x_{_{1}})=
      a_{_{k}}\,x_{_{0}}\,x_{_{1}}\,\bar{a_{_{k}}}=x_{_{0}}\,x_{_{k}}\,.$$
\item It follows from the braid relations and the case 2 that 
      $$b(x_{_{k}})=b\,a_{_{k}}\,b(x_{_{0}})=a_{_{k}}\,b\,
      a_{_{k}}(x_{_{0}})=a_{_{k}}\,b(x_{_{0}})=x_{_{k}}\,,$$
      and we get also $\,\bar{b}(x_{_{k}})\!=\!x_{_{k}}\,$.
\item For $\,k\!\not =\! 2$, one has $\,b_{_{1}}(x_{_{k}})\!=\!
      \bar{b_{_{1}}}(x_{_{k}})\!=\!x_{_{k}}\,$ by the braid relations. 
      When $\,k\!=\!2$, we get $\,b_{_{1}}(x_{_{2}})\!=\!x_{_{3}}\,$ and
      $\,\bar{b_{_{1}}}(x_{_{2}})\!=\!x_{_{2}}\,\bar{x_{_{3}}}\,x_{_{2}}\,$:
      $$\begin{array}{rcll}
      x_{_{2}}\,\bar{x_{_{3}}}\,x_{_{2}} & = & a_{_{2}}\,x_{_{1}}\,
      \bar{a_{_{2}}}\,b_{_{1}}\,a_{_{2}}\,\bar{x_{_{1}}}\,\bar{a_{_{2}}}\,
      \bar{b_{_{1}}}\,a_{_{2}}\,x_{_{1}}\,\bar{a_{_{2}}} & \\
      
      &=& a_{_{2}}\,x_{_{1}}\,b_{_{1}}\,a_{_{2}}\,\bar{b_{_{1}}}\,
      \bar{x_{_{1}}}\,b_{_{1}}\,\bar{a_{_{2}}}\,\bar{b_{_{1}}}\,x_{_{1}}\,
      \bar{a_{_{2}}} & \hbox{by {\it (T)}} \\
      
      &=& a_{_{2}}\,b_{_{1}}\,x_{_{1}}\,\bar{x_{_{2}}}\,x_{_{1}}\,
      \bar{b_{_{1}}}\,\bar{a_{_{2}}} & \hbox{by case 2} \\
      
      &=& a_{_{2}}\,b_{_{1}}\,\bar{a_{_{2}}}\,x_{_{1}}\,a_{_{2}}
      \bar{b_{_{1}}}\,\bar{a_{_{2}}} & \hbox{by case 2} \\
      
      &=& \bar{b_{_{1}}}\,a_{_{2}}\,b_{_{1}}\,x_{_{1}}\,\bar{b_{_{1}}}\,
      \bar{a_{_{2}}}\,b_{_{1}} & \hbox{by {\it (T)}} \\
      
      &=& \bar{b_{_{1}}}(x_{_{2}}) & \hbox{by case 2}\,.
      \end{array}$$
\item Each $\,y\!\in\!\{b_{_{2}},\ldots,b_{_{g-1}},c_{_{2,4}},\ldots,
      c_{_{2g-4,2g-2}},c_{_{1,2}}\}\,$ commutes with $\,x_{_{k}}\,$ for
      $\,k\!=\!2,2g,\ldots, 2g+n-3\,$ by 
      the braid relations. Therefore,  we get
      $\,y(x_{_{k}})\!=\!\bar{y}(x_{_{k}})\!=\!x_{_{k}}\,$.
\item Let $\,i\!\in\!\{2,2g,\ldots,2g+n-3\}$. Suppose 
      first that $\,i\!\geq\!k$. Then, if $\,m_{_{k}}\!=\!\bar{x_{_{1}}}
      (a_{_{k}})$, we have
      $$a_{_{i}}(x_{_{k}})=a_{_{i}}\,a_{_{k}}\,x_{_{1}}\,\bar{a_{_{k}}}\,
      \bar{a_{_{i}}}=a_{_{i}}\,x_{_{1}}\,m_{_{k}}\,
      \bar{a_{_{i}}}\,\bar{a_{_{k}}}\,.$$
      By the braid relations, one has
      $$m_{_{k}}=b\,\bar{a_{_{1}}}\,a_{_{2g+n-2}}\,\bar{b}(a_{_{k}})=b\,
      \bar{a_{_{1}}}\,a_{_{2g+n-2}}\,a_{_{k}}(b)=b\,a_{_{2g+n-2}}\,a_{_{k}}\,
      b(a_{_{1}})$$
      and the lantern relation $\,(L_{_{2g+n-2,1,k}})\,$ says that
      $$a_{_{2g+n-2}}\,c_{_{2g+n-2,1}}\,c_{_{1,k}}\,a_{_{k}}=c_{_{2g+n-2,k}}\,
      a_{_{1}}\,Y\,a_{_{1}}\,\bar{Y}$$
      where $\,Y=b\,a_{_{2g+n-2}}\,a_{_{k}}\,b$. Thus, we get
      $$m_{_{k}}=Y(a_{_{1}})=\bar{a_{_{1}}}\,\bar{c_{_{2g+n-2,k}}}\,
      a_{_{2g+n-2}}\,c_{_{2g+n-2,1}}\,c_{_{1,k}}\,a_{_{k}}\,,$$
      which implies by the braid relations $\,m_{_{k}}a_{_{i}}\!=\!a_{_{i}}
      m_{_{k}}\,$ since $\,i\!\geq\!k$. From this, one obtains
      $$a_{_{i}}(x_{_{k}})=a_{_{i}}\,x_{_{1}}\,\bar{a_{_{i}}}\,m_{_{k}}\,
      \bar{a_{_{k}}}=a_{_{i}}\,x_{_{1}}\,\bar{a_{_{i}}}\,\bar{x_{_{1}}}\,
      a_{_{k}}\,x_{_{1}}\,\bar{a_{_{k}}}=x_{_{i}}\,\bar{x_{_{1}}}\,x_{_{k}}\,.$$
      In particular, we have  $\,x_{_{k}}\!=\!x_{_{1}}\,\bar{x_{_{i}}}\,
      a_{_{i}}\,x_{_{k}}\,\bar{a_{_{i}}}\,$ and so:
      $$\begin{array}{rcll}
      \bar{a_{_{i}}}(x_{_{k}}) & = & \bar{a_{_{i}}}\,x_{_{1}}\,\bar{x_{_{i}}}\,
      a_{_{i}}\,x_{_{k}}\,\bar{a_{_{i}}}\,a_{_{i}} & \\

      &=& \bar{a_{_{i}}}\,x_{_{1}}\,a_{_{i}}\,
      \bar{a_{_{i}}}\,\bar{x_{_{i}}}\,a_{_{i}}\,x_{_{k}} & \\

      &=& x_{_{1}}\,\bar{x_{_{i}}}\,x_{_{1}}\,
      \bar{x_{_{1}}}\,x_{_{k}} & \hbox{by case 2} \\

      &=& x_{_{1}}\,\bar{x_{_{i}}}\,x_{_{k}}\,. &
      \end{array}$$
      Conclusion: $\,\left\{ \begin{array}{ccl}
      a_{_{i}}(x_{_{k}})=x_{_{i}}\,\bar{x_{_{1}}}\,x_{_{k}}\,, &
      \bar{a_{_{i}}}(x_{_{k}})=x_{_{1}}\,\bar{x_{_{i}}}\,x_{_{k}} &
      \hbox{if } \,i\geq k, \\
      a_{_{i}}(x_{_{k}})=x_{_{k}}\,\bar{x_{_{1}}}\,x_{_{i}} \,, &
      \bar{a_{_{i}}}(x_{_{k}})= x_{_{1}}\,\bar{x_{_{k}}}\,x_{_{i}} &
      \hbox{if } \,i\leq k.
      \end{array}\right.$

\end{list}

\vskip5mm\noindent
$\ast$ \underline{Case 4: $\,k\!=\!3$}.

\vskip3mm
\begin{list}{--}{\leftmargin7mm\itemsep3mm}
\item By the braid relations and the preceeding cases, we have:
      $$a_{_{1}}(x_{_{3}})=b_{_{1}}\,a_{_{1}}(x_{_{2}})=b_{_{1}}\,
      x_{_{2}}\,\bar{x_{_{0}}}\,\bar{b_{_{1}}}=x_{_{3}}\,\bar{x_{_{0}}}\,,$$
      $$\bar{a_{_{1}}}(x_{_{3}})=b_{_{1}}\,\bar{a_{_{1}}}(x_{_{2}})=b_{_{1}}\,
      x_{_{2}}\,x_{_{0}}\,\bar{b_{_{1}}}=x_{_{3}}\,x_{_{0}}\,,$$
      $$a_{_{2g+n-2}}(x_{_{3}})=b_{_{1}}\,a_{_{2g+n-2}}(x_{_{2}})=b_{_{1}}\,
      \bar{x_{_{0}}}\,x_{_{2}}\,\bar{b_{_{1}}}=\bar{x_{_{0}}}\,x_{_{3}}\,,$$
      $$\bar{a_{_{2g+n-2}}}(x_{_{3}})=b_{_{1}}\,\bar{a_{_{2g+n-2}}}(x_{_{2}})=
      b_{_{1}}\,x_{_{0}}\,x_{_{2}}\,\bar{b_{_{1}}}=x_{_{0}}\,x_{_{3}}\,.$$
\item The relations $\,${\it (T)}$\,$ and the case 3 prove that
      $$b(x_{_{3}})=b\,b_{_{1}}(x_{_{2}})=b_{_{1}}(x_{_{2}})=x_{_{3}}=
      \bar{b}(x_{_{3}}),$$
      and
      $$a_{_{2}}(x_{_{3}})=a_{_{2}}\,b_{_{1}}\,a_{_{2}}(x_{_{1}})=
      b_{_{1}}\,a_{_{2}}\,b_{_{1}}(x_{_{1}})=b_{_{1}}\,a_{_{2}}(x_{_{1}})
      =x_{_{3}}=\bar{a_{_{2}}}(x_{_{3}}).$$
\item One has $\,\bar{b_{_{1}}}(x_{_{3}})\!=\!x_{_{2}}\,$. On the 
      other hand, we get
      $$\begin{array}{rcll}
      b_{_{1}}(x_{_{3}}) & = & b_{_{1}}\,x_{_{2}}\,\bar{b_{_{1}}}\,
      \bar{x_{_{2}}}\,b_{_{1}}\,x_{_{2}}\,\bar{b_{_{1}}} & \hbox{by case 3} \\
      &=& x_{_{3}}\,\bar{x_{_{2}}}\,x_{_{3}}\,. & \\
      \end{array}$$
\item Using the braid relations and the case 3, we get $\,\bar{c_{_{2,4}}}
      (x_{_{3}})\!=\!x_{_{3}}\,\bar{x_{_{4}}}\,x_{_{3}}\,$:
      $$\begin{array}{rcl}
      x_{_{3}}\,\bar{x_{_{4}}}\,x_{_{3}} & = & b_{_{1}}\,x_{_{2}}\,
      \bar{b_{_{1}}}\,c_{_{2,4}}\,b_{_{1}}\,\bar{x_{_{2}}}\,\bar{b_{_{1}}}\,
      \bar{c_{_{2,4}}}\,b_{_{1}}\,x_{_{2}}\,\bar{b_{_{1}}} \\
      
      &=& b_{_{1}}\,x_{_{2}}\,c_{_{2,4}}\,b_{_{1}}\,\bar{c_{_{2,4}}}\,
      \bar{x_{_{2}}}\,c_{_{2,4}}\,\bar{b_{_{1}}}\,\bar{c_{_{2,4}}}\,x_{_{2}}\,
      \bar{b_{_{1}}} \\
      
      &=& b_{_{1}}\,c_{_{2,4}}\,x_{_{2}}\,\bar{x_{_{3}}}\,x_{_{2}}\,
      \bar{c_{_{2,4}}}\,\bar{b_{_{1}}} \\
      
      &=& b_{_{1}}\,c_{_{2,4}}\,\bar{b_{_{1}}}\,x_{_{2}}\,b_{_{1}}\,
      \bar{c_{_{2,4}}}\,\bar{b_{_{1}}} \\
      
      &=& \bar{c_{_{2,4}}}\,b_{_{1}}\,c_{_{2,4}}\,x_{_{2}}\,\bar{c_{_{2,4}}}\,
      \bar{b_{_{1}}}\,c_{_{2,4}} \\
      
      &=& \bar{c_{_{2,4}}}(x_{_{3}}).
      \end{array}$$
      On the other hand, we have 
      $\,c_{_{2,4}}(x_{_{3}})\!=\!x_{_{4}}\,$.
\item The braid relations assure that $\,y(x_{_{3}})\!\!=\!\!\bar{y}(x_{_{3}})
      \!\!=\!\!
      x_{_{3}}\,$ for all \linebreak[4]$\,y\!\in\!\{b_{_{2}},\ldots,b_{_{g-1}},
      c_{_{4,6}},\ldots,c_{_{2g-4,2g-2}}\}$.
\item For each $\,i\!\in\!\{2g,\ldots,2g+n-3\}$, one has by the case 3
      $$a_{_{i}}(x_{_{3}})=b_{_{1}}\,a_{_{i}}(x_{_{2}})=
      b_{_{1}}\,x_{_{i}}\,\bar{x_{_{1}}}\,
      x_{_{2}}\,\bar{b_{_{1}}}=x_{_{i}}\,\,\bar{x_{_{1}}}\,x_{_{3}}\,$$
      and
      $$\bar{a_{_{i}}}(x_{_{3}})=b_{_{1}}\,\bar{a_{_{i}}}(x_{_{2}})=
      b_{_{1}}\,x_{_{1}}\,\bar{x_{_{i}}}\,x_{_{2}}\,\bar{b_{_{1}}}=
      x_{_{1}}\,\bar{x_{_{i}}}\,x_{_{3}}\,.$$
\item Finally, we shall prove that $\,c_{_{1,2}}(x_{_{3}})\!=\!x_{_{3}}\,
      \bar{x_{_{2}}}\,x_{_{1}}\,\bar{x_{_{0}}}\,d_{_{n}}\,$.
      
      \vskip3mm\noindent
      The lantern relation $\,(L_{_{2g+n-2,1,2}})\,$ says
      $$a_{_{2g+n-2}}\,c_{_{2g+n-2,1}}\,c_{_{1,2}}\,a_{_{2}}=c_{_{2g+n-2,2}}\,
      \bar{X}\,a_{_{1}}\,X\,a_{_{1}}=c_{_{2g+n-2,2}}\,a_{_{1}}\,X\,a_{_{1}}
      \,\bar{X}$$
      where $\,X\!=\!b\,a_{_{2}}\,a_{_{2g+n-2}}\,b$, that is to say 
      $\,$($d_{_{n}}\!=\!c_{_{2g+n-2,1}}$):
      $$a_{_{2g+n-2}}\,c_{_{1,2}}\,\bar{a_{_{1}}}=c_{_{2g+n-2,2}}\,
      \bar{a_{_{2}}}\,\bar{d_{_{n}}}\,\bar{X}\,a_{_{1}}\,X\ \ \ 
      (\star)$$
      and
      $$c_{_{2g+n-2,2}}\,\bar{c_{_{1,2}}}=X\,\bar{a_{_{1}}}\,\bar{X}\,
      \bar{a_{_{1}}}\,a_{_{2}}\,d_{_{n}}\,a_{_{2g+n-2}}\ \ \ 
      (\star\star).$$
      Then, one can compute
      $$\begin{array}{rcll}
      \bar{x_{_{3}}}(c_{_{1,2}}) & = & b_{_{1}}\,a_{_{2}}\,b\,a_{_{2g+n-2}}\,
      \bar{a_{_{1}}}\,\bar{b}\,\bar{a_{_{2}}}\,\bar{b_{_{1}}}(c_{_{1,2}}) & \\
      
      &=& b_{_{1}}\,a_{_{2}}\,b\,a_{_{2g+n-2}}\,c_{_{1,2}}\,\bar{a_{_{1}}}\,
      \bar{b}\,\bar{a_{_{2}}}(b_{_{1}}) & \hbox{by {\it (T)}} \\
      
      &=& b_{_{1}}\,a_{_{2}}\,b\,c_{_{2g+n-2,2}}\,\bar{a_{_{2}}}\,
      \bar{d_{_{n}}}\,\bar{X}\,a_{_{1}}\,X\,\bar{b}\,
      \bar{a_{_{2}}}(b_{_{1}}) & \hbox{by }\,(\star) \\
      
      &=& b_{_{1}}\,a_{_{2}}\,b\,c_{_{2g+n-2,2}}\,\bar{a_{_{2}}}\,
      \bar{d_{_{n}}}\,\bar{X}\,a_{_{1}}\,b\,a_{_{2g+n-2}}(b_{_{1}}) &  \\
      
      &=& b_{_{1}}\,c_{_{2g+n-2,2}}\,\bar{b}\,a_{_{2}}\,b\,\bar{X}(b_{_{1}})
       & \hbox{by {\it (T)}}  \\
      
      &=& b_{_{1}}\,c_{_{2g+n-2,2}}\,\bar{b}\,a_{_{2}}\,b\,\bar{b}\,
      \bar{a_{_{2}}}\,\bar{a_{_{2g+n-2}}}\,\bar{b}(b_{_{1}}) &  \\
            
      &=& b_{_{1}}\,\bar{b_{_{1}}}(c_{_{2g+n-2,2}}) & \hbox{by {\it (T)}} \\
            
      &=& c_{_{2g+n-2,2}}\,. &
      \end{array}$$
      Thus, we get 
      $$\begin{array}{rcll}
      c_{_{1,2}}(x_{_{3}}) &= & c_{_{1,2}}\,x_{_{3}}\,\bar{c_{_{1,2}}} & \\
     
      &=& x_{_{3}}\,\bar{x_{_{3}}}\,c_{_{1,2}}\,x_{_{3}}\,\bar{c_{_{1,2}}} & \\
     
      &=& x_{_{3}}\,c_{_{2g+n-2,2}}\,\bar{c_{_{1,2}}} & \\
      
      &=& x_{_{3}}\,X\,\bar{a_{_{1}}}\,\bar{X}\,\bar{a_{_{1}}}\,a_{_{2}}\,
      a_{_{2g+n-2}}\,d_{_{n}} & \hbox{by }\,(\star\star) \\
      
      &=& x_{_{3}}\,b\,a_{_{2}}\,a_{_{2g+n-2}}\,b\,\bar{a_{_{1}}}\,\bar{b}\,
      \bar{a_{_{2}}}\,\bar{a_{_{2g+n-2}}}\,\bar{b}\,
      a_{_{2}}\,\bar{x_{_{0}}}\,d_{_{n}} &  \\
      
      &=& x_{_{3}}\,b\,a_{_{2g+n-2}}\,a_{_{2}}\,\bar{a_{_{1}}}\,\bar{b}\,
      a_{_{1}}\,\bar{a_{_{2}}}\,\bar{a_{_{2g+n-2}}}\,\bar{b}\,
      a_{_{2}}\,\bar{x_{_{0}}}\,d_{_{n}} &  \\
      
      &=& x_{_{3}}\,b\,\bar{x_{_{0}}}\,\bar{b}\,
      \bar{a_{_{2}}}\,b\,x_{_{0}}\,\bar{b}\,
      a_{_{2}}\,\bar{x_{_{0}}}\,d_{_{n}} & \hbox{by {\it (T)}} \\
      
      &=& x_{_{3}}\,\bar{x_{_{1}}}\,\bar{a_{_{2}}}\,x_{_{1}}\,a_{_{2}}\,
      \bar{x_{_{0}}}\,d_{_{n}} & \\
      
      &=& x_{_{3}}\,\bar{x_{_{1}}}\,x_{_{1}}\,\bar{x_{_{2}}}\,x_{_{1}}\,
      \bar{x_{_{0}}}\,d_{_{n}} & \hbox{by case 2}\\
 
      &=& x_{_{3}}\,\bar{x_{_{2}}}\,x_{_{1}}\,
      \bar{x_{_{0}}}\,d_{_{n}}\,. &
      \end{array}$$
      It follows from this that
      $$\bar{c_{_{1,2}}}(x_{_{3}})\!=\bar{c_{_{1,2}}}\,c_{_{1,2}}\,x_{_{3}}\,
      \bar{c_{_{1,2}}}\,\bar{d_{_{n}}}\,x_{_{0}}\,\bar{x_{_{1}}}\,x_{_{2}}\,
      c_{_{1,2}}=x_{_{3}}\,\bar{d_{_{n}}}\,x_{_{0}}\,\bar{x_{_{1}}}\,x_{_{2}}\,
      .$$

\end{list}

\vskip5mm\noindent
$\ast$ \underline{Case 5: $\,k\!\in\!\{4,5,\ldots,2g-1\}$}.

\vskip3mm\noindent
In order to simplify the notation, let us denote
$$e_{_{3}}=b_{_{1}}\,,\ \,e_{_{4}}=c_{_{2,4}}\,,\ \,e_{_{5}}=b_{_{2}}\,,\ \,
\ldots\,,\ \,e_{_{2g-2}}=c_{_{2g-4,2g-2}}\,,\ \ e_{_{2g-1}}=b_{_{g-1}}\,,$$
so that, for $\,i\!\in\!\{3,\ldots,2g-1\}$, 
$\,x_{_{i}}\!=\!e_{_{i}}(x_{_{i-1}})$.

\vskip3mm
\begin{list}{--}{\leftmargin7mm\itemsep3mm}
\item Then, one has by the braid relations and the case 4:
      $$a_{_{1}}(x_{_{k}})=e_{_{k}}\,e_{_{k-1}}\cdots 
      e_{_{4}}\,a_{_{1}}(x_{_{3}})=e_{_{k}}\cdots e_{_{4}} 
      x_{_{3}}\,\bar{x_{_{0}}}\,\bar{e_{_{4}}}\cdots\bar{e_{_{k}}}=x_{_{k}}\,
      \bar{x_{_{0}}}\,.$$
      Likewise, we get
      $$\,\bar{a_{_{1}}}(x_{_{k}})=x_{_{k}}\,x_{_{0}}\,,\ \ \ 
      a_{_{2g+n-2}}(x_{_{k}})=\bar{x_{_{0}}}\,x_{_{k}}\,,\ \ \ 
      \bar{a_{_{2g+n-2}}}(x_{_{k}})=x_{_{0}}\,x_{_{k}}\,,$$ 
      $$\hbox{and }\ \,b(x_{_{k}})=\bar{b}(x_{_{k}})=x_{_{k}}=
      a_{_{2}}(x_{_{k}})=\bar{a_{_{2}}}(x_{_{k}})\,.$$
\item For $\,i\!\in\!\{3,4,\ldots,2g-1\},\ \,i\!<\!k$, one obtains, 
      using the braid relations, $\,e_{_{i}}(x_{_{k}})\!=\bar{e_{_{i}}}(x_{_{k}})=
      x_{_{k}}\,$:

      \begin{center}
      $e_{_{i}}(x_{_{k}})$ \parbox[t]{98mm}{$=e_{_{k}}\cdots e_{_{i}}\,
      e_{_{i+1}}\,e_{_{i}}\cdots e_{_{3}}(x_{_{2}})=e_{_{k}}\cdots e_{_{i+1}}\,
      e_{_{i}}\,e_{_{i+1}}\cdots e_{_{3}}(x_{_{2}})$
      $=e_{_{k}}\cdots e_{_{3}}(x_{_{2}})=x_{_{k}}\,.$}
      \end{center}

      \noindent
      For $\,i\!>\!k+1$, $\,e_{_{i}}\,$ commutes with 
      $\,e_{_{k}}\,,\ldots,\,e_{_{4}}\,$ and $\,x_{_{3}}\,$, thus we also have
      $$e_{_{i}}(x_{_{k}})=\bar{e_{_{i}}}(x_{_{k}})=x_{_{k}}\,\ 
      (i>k+1)\ \ \ (\ast).$$
\item One has $\,e_{_{k+1}}(x_{_{k}})\!=\!x_{_{k+1}}\,$. Let us prove by
      induction on $k$ that $\,\bar{e_{_{k+1}}}(x_{_{k}})\!=\!x_{_{k}}\,
      \bar{x_{_{k+1}}}\,x_{_{k}}\,$. We have seen in case 4 that this equaliy
      is satisfied at the rank $\,k\!=\!3$. Suppose it is true at the rank 
      $k\!-\!1$,  $\,4\!\leq\!k\leq 2g-2$. Then, we get:
      $$\begin{array}{rcl}
      x_{_{k}}\,\bar{x_{_{k+1}}}\,x_{_{k}} & = & e_{_{k}}\,x_{_{k-1}}\,
      \bar{e_{_{k}}}\,e_{_{k+1}}\,e_{_{k}}\,\bar{x_{_{k-1}}}\,\bar{e_{_{k}}}\,
      \bar{e_{_{k+1}}}\,e_{_{k}}\,x_{_{k-1}}\,\bar{e_{_{k}}} \\
      
      &=& e_{_{k}}\,x_{_{k-1}}\,e_{_{k+1}}\,e_{_{k}}\,\bar{e_{_{k+1}}}\,
      \bar{x_{_{k-1}}}\,e_{_{k+1}}\,\bar{e_{_{k}}}\,\bar{e_{_{k+1}}}\,
      x_{_{k-1}}\,\bar{e_{_{k}}} \ \  \hbox{by {\it (T)}} \\

      &=& e_{_{k}}\,e_{_{k+1}}\,x_{_{k-1}}\,e_{_{k}}\,\bar{x_{_{k-1}}}\,
      \bar{e_{_{k}}}\,x_{_{k-1}}\,\bar{e_{_{k+1}}}\,
      \bar{e_{_{k}}}\ \ \ \  \hbox{by }\ (\ast) \\

      &=& e_{_{k}}\,e_{_{k+1}}\,x_{_{k-1}}\,\bar{x_{_{k}}}\,
      x_{_{k-1}}\,\bar{e_{_{k+1}}}\,\bar{e_{_{k}}}  \\

      &=& e_{_{k}}\,e_{_{k+1}}\,\bar{e_{_{k}}}\,x_{_{k-1}}\,e_{_{k}}\,
      \bar{e_{_{k+1}}}\,\bar{e_{_{k}}} \ \ \ \hbox{by inductive 
      hypothesis} \\

      &=& \bar{e_{_{k+1}}}\,e_{_{k}}\,e_{_{k+1}}\,x_{_{k-1}}\,\bar{e_{_{k+1}}}\,
      \bar{e_{_{k}}}\,e_{_{k+1}} \ \ \ \  \hbox{by {\it (T)}} \\

      &=& \bar{e_{_{k+1}}}\,e_{_{k}}\,x_{_{k-1}}\,\bar{e_{_{k}}}\,e_{_{k+1}}
      \ \ \ \  \hbox{by }\ (\ast) \\

      &=& \bar{e_{_{k+1}}}(x_{_{k}}).
      \end{array}$$
\item This last relation implies $\,x_{_{k}}\!=\!x_{_{k-1}}\,\bar{e_{_{k}}}\,
      \bar{x_{_{k-1}}}\,e_{_{k}}\,x_{_{k-1}}\,$. Thus, we get
      $$e_{_{k}}(x_{_{k}})=e_{_{k}}\,x_{_{k-1}}\,\bar{e_{_{k}}}\,
      \bar{x_{_{k-1}}}\,e_{_{k}}\,x_{_{k-1}}\,\bar{e_{_{k}}}=x_{_{k}}\,
      \bar{x_{_{k-1}}}\,x_{_{k}}\,.$$
      On the other hand, one has $\,\bar{e_{_{k}}}(x_{_{k}})\!=\!x_{_{k-1}}\,$.
\item For $\,i\!\in\!\{2g,\ldots,2g+n-3\}$, we have, by the braid
      relations and the cases 2, 3 and 4:
      $$a_{_{i}}(x_{_{k}})=e_{_{k}}\cdots e_{_{4}}\,a_{_{i}}(x_{_{3}})=
      e_{_{k}}\cdots e_{_{4}}\,x_{_{i}}\,\bar{x_{_{1}}}\,x_{_{3}}\,
      \bar{e_{_{4}}}\cdots\bar{e_{_{k}}}=x_{_{i}}\,\bar{x_{_{1}}}\,x_{_{k}}\,,$$
      and likewise, we get $\,\bar{a_{_{i}}}(x_{_{k}})\!=x_{_{1}}\,
      \bar{x_{_{i}}}\,x_{_{k}}\,$.
\item Finally, since $\,c_{_{1,2}}(x_{_{3}})\!=\!x_{_{3}}\,\bar{x_{_{2}}}\,
      x_{_{1}}\,\bar{x_{_{0}}}\,d_{_{n}}$, it follows from the braid
      relations and the preceeding cases that $\,c_{_{1,2}}(x_{_{k}})\!=\!x_{_{k}}\,
      \bar{x_{_{2}}}\,x_{_{1}}\,\bar{x_{_{0}}}\,d_{_{n}}$. In the 
      same way, we get $\,\bar{c_{_{1,2}}}(x_{_{k}})\!=\!x_{_{k}}\,
      \bar{d_{_{n}}}\,x_{_{0}}\,\bar{x_{_{1}}}\,x_{_{2}}$.
\end{list}
\eproof

\vskip3mm\noindent
{\bf Proof of proposition~\ref{h}.\ \,} If $\,\pi:\mathbf{\bf 
Z}\times \pi_{1}(\Sigma_{g,n-1},p)\rightarrow 
\pi_{1}(\Sigma_{g,n-1},p)\,$ denotes the projection, the loops 
$\,\pi\circ h_{g,n}(x_{_{0}}),\ldots,\pi\circ 
h_{g,n}(x_{_{2g+n-3}})\,$ form a basis of the free group 
$\,\pi_{1}(\Sigma_{g,n-1},p)$. Thus, $F$, the subgroup of $\,\ker g_{2}\,$ 
generated by $\,x_{_{0}},\ldots,x_{_{2g+n-3}}\,$ is free of rank 
$2g+n-2$ and the restriction of $\,\pi\circ h_{g,n}\,$ to this subgroup is an
isomorphism.

\noindent
Now, for all element $x$ of $\,\ker g_{2}$, there are by 
lemma~\ref{gen-ker g2} an integer $k$ and an element $f$ of $F$ such 
that $\,x\!=\!d_{_{n}}^k\,f\,$ ($d_{_{n}}\,$ is central in $\,\ker 
g_{2}\,$). Then, one has $\,h_{g,n}(x)\!=\!\bigr(k,\pi\circ 
h_{g,n}(x)\bigl)\,$ and therefore, $\,h_{g,n}\,$ is one to one. But $h_{g,n}\,$
is also onto. This concludes the proof.

\eproof

\noindent
{\bf Proof of theorem~\ref{principaltheorem}.\ \,}
In section~\ref{negalun}, we proved that $\Phi_{g,1}\,$ is an isomorphism.
Thus, by the five-lemma, proposition~\ref{h} and an inductive argument, $\Phi_{g,n}\,$ is 
an isomorphism for all $\,n\geq 1$. In order to conclude the proof, 
it remains to look at the case $\,n\!=\!0$.

Since all spin maps are conjugate in $\,\mathcal{M}_{g,1}$, $\,\ker f_{2}\,$ is 
normally\linebreak[4] generated by $\,\tau_{\delta_{1}}\,$ and 
$\,\tau_{\alpha_{1}}\tau_{\alpha_{2g-1}}^{-1}$. Thus, considering once more 
the commutative diagram

\diagram[size=2.5em]
 1 & \rto &   \ker g_{2}   & \rto & G_{g,1} &  \rTo^{g_{2}} && G_{g,0} & \rto & 1  \\
   &      & \dto_{h_{g,1}}  &   & \dto_{\Phi_{g,1}}^{\approx} & && 
   \dto_{\Phi_{g,0}} & &\\
 1 & \rto & \mathbf{Z}\times\pi_{1}(\Sigma_{g,0},p)  & \rTo^{f_{1}} &  \mathcal{M}_{g,1} & 
 \rto^{f_{2}} && \mathcal{M}_{g,0} & \rto & 1 \\
\enddiagram

\vskip3mm\noindent
and recalling that $\,\ker g_{2}\,$ is normally generated by 
$\,d_{_{1}}\,$ and $\,a_{_{1}}\,\bar{a_{_{2g-1}}}\,$ 
(lemma~\ref{normal}), we conclude 
that $\,h_{g,1}\,$ is still an isomorphism. So, we get that 
$\Phi_{g,0}\,$ is an isomorphism.

\eproof


\vskip7mm\noindent
{\bf Acknowledgement. } This paper originates from discussions I had with 
Catherine Labru\`ere. I want to thank her. 


\bibliographystyle{amsplain}

\end{document}